\documentclass[12pt,reqno]{amsart}
\usepackage{amssymb,amsfonts,amsthm,amsmath,mathrsfs}
\usepackage{cite}
\usepackage[shortlabels]{enumitem}
\usepackage[left=1 in,top=1 in,right=1 in, bottom=1 in]{geometry}
\usepackage{graphicx}
\usepackage{color}
\usepackage{geometry}
\usepackage[pagebackref=false]{hyperref}
\usepackage{tcolorbox}

\def\d{{\rm d}}
\def\eqdef{\stackrel{\rm def}{=}}
\definecolor{darkred}{rgb}{.70,.12,.20}

\definecolor{darkgreen}{rgb}{.20,.52,.14}

\definecolor{byz}{rgb}{.44,.16,.39}

\numberwithin{equation}{section}

\newtheorem{theorem}{Theorem}[section]
\newtheorem{lemma}[theorem]{Lemma}
\newtheorem{definition}[theorem]{Definition}
\newtheorem{corollary}[theorem]{Corollary}
\newtheorem{assumption}[theorem]{Assumption}

\theoremstyle{remark}
\newtheorem{remark}[theorem]{Remark}
\newtheorem{example}[theorem]{\bf{Example}}

\usepackage{todonotes}

\newcommand{\osc}{\mathop{\mathrm{osc}}}

\newcommand{\varep}{\varepsilon}

\newcommand{\beq}{\begin{equation}}
\newcommand{\eeq}{\end{equation}}
\newcommand{\beqs}{\begin{equation*}}
\newcommand{\eeqs}{\end{equation*}}
\newcommand{\ba}{\begin{array}}
\newcommand{\ea}{\end{array}}
\newcommand{\beas}{\begin{eqnarray*}}
\newcommand{\eeas}{\end{eqnarray*}}
\newcommand{\bea}{\begin{eqnarray}}
\newcommand{\eea}{\end{eqnarray}}
\newcommand{\bal}{\begin{align}}
\newcommand{\eal}{\end{align}}

\newcommand{\bals}{\begin{align*}}
\newcommand{\eals}{\end{align*}}

\newcommand{\tnum}{\rm(\roman*)}
\newcommand{\rnum}{\rm(\alph*)}

\newcommand{\R}{\ensuremath{\mathbb R}}

\newcommand{\bds}{\begin{displaystyle}}
\newcommand{\eds}{\end{displaystyle}}

\newcommand{\remove}[1]{} 
\renewcommand{\remove}[1]{#1} 

\definecolor{darkred}{rgb}{.70,.12,.20}

\definecolor{darkgreen}{rgb}{.20,.52,.14}

\title{Linear non-divergence elliptic equations in a bounded, infinitely winding planar domain}

\author{Luan Hoang$^{1,*}$}
\address{$^{1}$Department of Mathematics and Statistics,
Texas Tech University,
1108 Memorial Circle, Lubbock, TX 79409--1042, U. S. A.}
\email{luan.hoang@ttu.edu}

\author{Akif Ibragimov$^{1,2}$}
\address{$^{2}$Oil and Gas Research Institute, Russian Academy of Science,  3 Gubkin Street, Moscow, 119333, Russia}
\email{ilya1sergey@gmail.com}

\thanks{$^*$Corresponding author.}

\date{\today}



\begin{document}

\begin{abstract} 
We study the second order elliptic equations of non-divergence form in a planar domain with complicated geometry. In this case the domain winds around a fixed circle infinitely many times and converges to it when the rotating angle goes to infinity. For the homogeneous equation and the homogeneous Dirichlet boundary condition,  in the case of bounded drifts, we prove that the maximum of the solution on the cross-section corresponding to a given rotating angle either grows  or decays exponentially as the angle goes to infinity. Results for the oscillation and its asymptotic estimates are also  obtained for inhomogeneous Dirichlet data. If the drift is unbounded but does not grow to infinity too fast, then the above maximum also goes to either  zero or  infinity. For the inhomogeneous equation, we obtain the estimates in the case of bounded forcing functions. Moreover, we establish the uniqueness of the solution and its continuous dependence on the boundary data and the forcing function. 
\end{abstract}

\maketitle 
\tableofcontents 

\pagestyle{myheadings}\markboth{\sc L.~Hoang and A.~Ibragimov}
{\sc Non-divergence elliptic equations in a bounded, infinitely winding planar domain}


\section{Introduction}\label{intro}
The aim of this paper is to study the non-divergence elliptic equations of the second order  in domains with complicated geometry. In the existing literature of the qualitative study of classical solutions in bounded domains, the solutions are usually characterized by the ``short" Euclidean distance to the boundary points.
Our domain's geometry, however, is structured by its ``length" which is infinitely ``long" even though the domain itself is bounded. Therefore, our problem is close, in spirit, to the study in unbounded domains like strips or cones with the limiting ``point" of interest being infinity. In the current work, we restrict ourselves to a two-dimensional domain and the set of limiting points is a closed curve where the solution's values are not specified. Moreover, we focus on the non-divergence,  elliptic equations. This situation requires a careful analysis which is the goal of this paper.

From the technical point of view, our essential tools will be the Maximum Principle and the Growth Lemma of Landis-type. Regarding the latter technique, in his papers \cite{Landis1963,Landis1968,Landis1968b}, see also \cite[Chapter 1]{LandisBook}, Landis used metric characteristics of the domain in terms of capacity.  This allowed him to connect the qualitative features of the solutions of elliptic and parabolic equations to the measure theory via  the corresponding relative capacity. In particular, this was used to study asymptotic properties of the solutions in unbounded domains. 
This type of lemma with capacity was also used  to prove the Wiener criteria for elliptic equations in divergence form \cite{Maz1967}, see also the book \cite{MazBook2018}.
In our version of the Growth Lemma, we use an explicit formula of the growth/contraction factor without involving the capacity concept. With this, we prove that the solution near the limiting set -- a circle in this case -- is stabilizing exponentially or its maximum on the cross-sections grows exponentially to infinity.
Moreover,  in  \cite{Landis1968}, the lower order terms either are ignored or subject to some structural  condition due to the unboundedness of the domain.
Here, we are able to deal with them -- in the form of drifts -- even unbounded ones under some constraints on their growth without any structural conditions. Regarding non-divergence equations with  unbounded drifts, other work such as \cite{LadyUral1988,Safonov2010,Safonov2018} have different settings and topics. For example, the paper \cite{Safonov2010} is focused on the Harnack inequality and Hopf--Olenik estimates, while our current work is focused on global and asymptotic estimates.

From the application point of view, this article takes up the elliptic problems from the diffusion-transport models for fluid flows in porous media \cite{HI3,HI4}. In those models, Einstein's probability method for the Brownian motion \cite{Einstein1905} was generalized to derive a nonlinear equation in non-divergence form. 
Here, we focus on its linearized equation which, as showed in \cite{HI3,HI4}, will play important roles in the analysis of the original nonlinear problem.

As far as the results are concerned, we will obtain global and asymptotic estimates for the sub-solutions/solutions of both homogeneous and inhomogeneous problems with the Dirichlet boundary data.  
For the former problem, both cases of bounded and unbounded drifts are studied. In a certain situation for the latter problem, additional oscillation estimates are also derived.
A sub-solution/solution's estimates are often in the form of a dichotomy: its maximum on the cross-sections either goes to zero or infinity as the rotating angle goes to infinity. We also establish the uniqueness of the solutions in a certain class of functions and their continuous dependence, globally and asymptotically,  on the boundary data and forcing functions. 

The paper is organized as follows.
In section \ref{domsec}, we describe the geometry of the domain. The main assumptions are \eqref{rnin} and \eqref{bassum1}.
Under additional Assumption \ref{bassum2}, we can characterize how far a point $x\in U$ goes along the domain, as it winds around a fixed circle, by introducing the ``arc-distance" ${\mathscr L}(x)$ in Definition \ref{arcdist}. Our results can be expressed with more geometric meaning using this distance ${\mathscr L}(x)$. The main relations between $\theta$ and ${\mathscr L}(x)$ are in Lemma \ref{lengths}.
In section \ref{glsec}, we obtain a couple of Growth Lemmas of Landis-type. The first is a general result in the $n$-dimensional space -- Lemma \ref{GGL0}. It is then applied to a planar domain in an annulus -- Lemma \ref{GLthick}. Since we do not require the domain to be narrow, the latter growth lemma will help us obtain global estimates for the sub-solutions and solutions.
In section \ref{estsec}, we estimate the solutions of the homogeneous problem. The dichotomy in the discrete form is presented in Lemma \ref{cyls}, which will  serve both cases of bounded and unbounded drifts.
The case of bounded drifts is treated in subsection \ref{bddrift} with Theorem \ref{thm1}, Corollaries \ref{cor1} and \ref{cor2}. Estimates are expressed in terms of the angle $\theta$ in Theorem \ref{thm1}, and arc-distance ${\mathscr L}(x)$ in Corollary \ref{cor1}.
The main statements are in the  continuum dichotomy form.
The case of unbounded drifts is treated in subsection \ref{ubdrift} with Theorem \ref{thm3}. Briefly speaking, if the drift does not grow too fast, see Assumption \ref{bmassum}, then we have the dichotomy: the maximum of the solution's absolute value on each cross-section corresponding to any given angle $\theta$ goes to zero or infinity as $\theta\to\infty$.
In section \ref{inhomsec}, we study the inhomogeneous problems in the case of bounded drifts. The case without a forcing function, we have estimates for the solutions in Theorem \ref{thmih1}, and estimates for their oscillations in Corollary \ref{corih1}. The case of non-zero forcing functions is treated in Theorem \ref{thmih3}. Note that the asymptotic estimates for the solutions only depend on the asymptotic behavior of the  boundary data and, when applicable, the forcing function.   
In section \ref{unicd}, we study the uniqueness and continuous dependence of the solution on the boundary data, see Theorem \ref{ThmDep1}, and, additionally, the forcing function, see Theorem \ref{ThmDep2}.
The ending Remark \ref{more} lists some possible future developments of this work. 

\section{Description of the domain}\label{domsec}

Our domain of interest, $U$, is an open, bounded subset of $\R^2$ with some special features that we describe below. Recall the standard change of variables for the polar coordinates 
\beqs
X(r,\theta)=(r\cos\theta,r\sin\theta)\text{ for }r\ge 0,\ \theta\in\R.
\eeqs

\begin{figure}[ht]
    \centering
    \includegraphics[scale=.35]{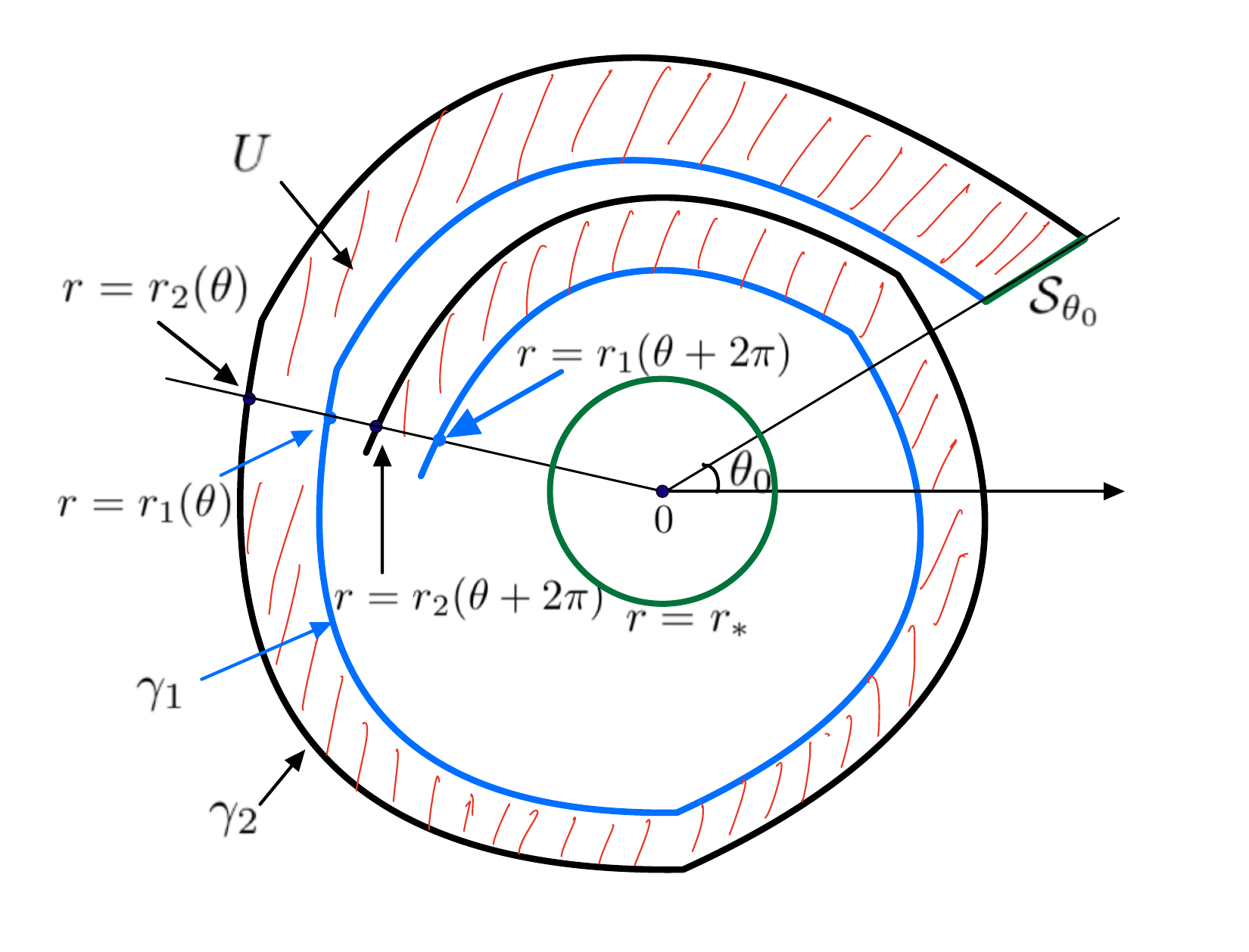}
    \caption{The domain $U$.}
    \label{PicDom}
 \end{figure}

We fix a number $\theta_0\in\R$ throughout the paper. Let $r_1$ and $r_2$ be positive, continuous functions on $[\theta_0,\infty)$ that satisfy 
\beq \label{rnin}
r_2(\theta+2\pi) <r_1(\theta)<r_2(\theta) \text{ for all } \theta\ge \theta_0.
\eeq
Define
\beqs
V=\{(r,\theta):\theta>\theta_0, r_1(\theta)<r<r_2(\theta)\},\quad U=X(V),
\eeqs
\beqs
\widetilde V=\{(r,\theta):\theta\ge \theta_0, r_0(\theta)\le r \le r_1(\theta)\}, \quad \widetilde U=X(\widetilde V),
\eeqs
see Fig. \ref{PicDom}. 
Thanks to \eqref{rnin}, one also has
\beq \label{rnin2}
r_1(\theta+2\pi) <r_2(\theta+2\pi)<r_1(\theta) \text{ for all } \theta\ge \theta_0.
\eeq
Therefore, the condition \eqref{rnin} guarantees that $X$ is a bijection from $U$ to $V$, and from $\widetilde U$ to $\widetilde V$.  In other words, $U$ and $\widetilde U$ each does not ``self-intersect". 
 
Define 
\beqs
\gamma=\gamma_1\cup\gamma_2, \text{ where }
\gamma_i=\{X(r,\theta):\theta\ge \theta_0, r=r_i(\theta)\} \text{ for }i=1,2.
\eeqs
For any number $\theta\ge \theta_0$, denote the line segment
\beqs
\mathcal S_\theta=\{ x=X(r,\theta): r_1(\theta)\le r\le r_2(\theta)\}.
\eeqs
This is the cross-section of $\widetilde U$ corresponding to the angle $\theta$.
In general, for any subset $I$ of $[\theta_0,\infty)$, denote
\beqs
U_I=\{X(r,\theta):\theta\in I,r_1(\theta)<r<r_2(\theta)\},\quad 
\widetilde U_I=\{X(r,\theta):\theta\in I, r_1(\theta)\le r\le r_2(\theta)\}.
\eeqs

Define the curve that is the convex combination of $\gamma_1$ and $\gamma_2$ with the ratio $\lambda\in[0,1]$ by
\beq \label{Xlam}
X_\lambda(\theta)=X(R_\lambda(\theta),\theta) \text{ for } \theta\in [\theta_0,\infty), \text{ where }
R_\lambda(\theta)=(1-\lambda)r_1(\theta)+ \lambda r_2(\theta).
\eeq
Clearly, $X_\lambda(\theta)\in \widetilde U$. Conversely, for $x\in \widetilde U$, there exists a unique pair $(r,\theta)\in \widetilde V$ such that $x=X(r,\theta)$. In particular, $r_1(\theta)\le r\le r_2(\theta)$, hence there is a unique number $\lambda\in[0,1]$ such that $r=R_\lambda(\theta)$. This implies  $x=X_\lambda(\theta)$. Therefore,
\beq\label{Xbi}
\text{the mapping $(\lambda,\theta)\mapsto X_\lambda(\theta)$ is a bijection from  $[0,1]\times[\theta_0,\infty)$ to $\widetilde U$.} 
\eeq

We assume throughout that the curves $\gamma_1$ and $\gamma_2$, roughly speaking,  converge to a fixed circle as $\theta\to\infty$. More precisely, 
there exists a number $r_*>0$ such that   
\beq \label{bassum1} 
\lim_{\theta\to\infty} r_1(\theta)=\lim_{\theta\to\infty} r_2(\theta)=r_*.
\eeq 

Because of the continuity of $r_1$, $r_2$ and assumption \eqref{bassum1}, there exists a number $\bar R>0$ such that
\beq\label{rrbound}
r_1(\theta)< r_2(\theta)\le \bar R\text{ for all }\theta\ge \theta_0.
\eeq
Denote the circle of radius $r_*$ centered at $0$ by ${\mathcal C}_*$.
By \eqref{rnin} and \eqref{rnin2}, for $i=1,2$ and $\theta\ge \theta_0$, the sequence $(r_i(\theta+2\pi k))_{k=0}^\infty$ is decreasing. Then we have 
\beqs  
r_*=\lim_{k\to\infty} r_i(\theta+2\pi k)=\inf \{r_i(\theta+2\pi k):\text{ for integers }k\ge 0\}.
\eeqs 
Hence, by \eqref{rnin}, 
\beqs
r_2(\theta)>r_1(\theta)>r_2(\theta+2\pi)\ge r_*\text{ for all }
\theta\ge \theta_0.
\eeqs 
Consequently, $\widetilde U\cap {\mathcal C}_*=\emptyset$.
Therefore, the closure of $U$, which is 
$\bar U=\widetilde U\cup {\mathcal C}_*$, is not the same as $\widetilde U$.

The main results below are obtained under the above assumption \eqref{rnin} and \eqref{bassum1}.
However, when the curves $\gamma_1$ and $\gamma_2$ have more regularity, more geometric information about the domain will be integrated into our analysis.
We consider such a case that gives rise to the notion of ``global length".

\begin{assumption}\label{bassum2}
The functions $r_1$ and $r_2$ are absolutely continuous on $[\theta_0,\infty)$, and there exists a number $\Theta\ge \theta_0$ such that  $r_1',r'_2\in L^\infty(\Theta,\infty)$.
\end{assumption}

A special case of Assumption \ref{bassum2} is when each norm $\|r_i'\|_{L^\infty(\theta,\infty)}$  converges to zero as $\theta\to\infty$. Thus, in this case, $\gamma_1$ and $\gamma_2$ converge to the circle $\mathcal C_*$, as $\theta\to\infty$, more strongly than the convergence in  \eqref{bassum1}.

Under Assumption \ref{bassum2}, for each $i=1,2$, thanks to the absolute continuity of $r_i$, the derivative $r_i'(\theta)$ exists for almost all $\theta\in[\theta_0,\infty)$ and 
\beq \label{rloc}
r_i'\in L^1_{\rm loc}([\theta_0,\infty)).
\eeq 
Define two numbers 
\begin{align}\label{muzero}
    \mu_0&=\max\{\|r_1'\|_{L^\infty(\Theta,\infty)},\|r_2'\|_{L^\infty(\Theta,\infty)}\}>0,\\
    \label{mustar}
     \mu_*&=\max\left\{\lim_{\theta\to\infty}\|r_1'\|_{L^\infty(\theta,\infty)},\lim_{\theta\to\infty}\|r_2'\|_{L^\infty(\theta,\infty)}\right\}\ge 0.
\end{align}
 For $\lambda\in[0,1]$ and $\theta\ge \theta_0$, denote the arc-length of the curve $X_\lambda([\theta_0,\theta])$ in \eqref{Xlam} by $s_\lambda(\theta)$, i.e.,
\beq\label{slform}
s_\lambda(\theta)=\int_{\theta_0}^\theta |X_\lambda'(\eta)| \d \eta, \text{ where the integral is the Lebesgue one.}
\eeq

\begin{definition}\label{arcdist}
Under Assumption \ref{bassum2}, for any $x\in \widetilde U$, define ${\mathscr L}(x)=s_\lambda(\theta)$, where $(\lambda,\theta)$ is the unique pair in $ [0,1]\times[\theta_0,\infty)$ such that $x=X_\lambda(\theta)$, see \eqref{Xbi}.
We call ${\mathscr L}(x)$ the arc-distance from $\mathcal S_{\theta_0}$ to $x$ along the domain $U$.
\end{definition}

Clearly, $x\mapsto{\mathscr L}(x)$ is a function from $\widetilde U$ to $[0,\infty)$.
We have the following comparisons between $s_\lambda(\theta)$ and ${\mathscr L}(x)$ with $\theta$. 

\begin{lemma}\label{lengths}
Under Assumption \ref{bassum2},  one has the following estimates with $x=X_\lambda(\theta)\in \widetilde U$. 
\begin{enumerate}[label=\tnum]
    \item For any $\lambda\in[0,1]$ and $\theta\ge \theta_0$,
    \beq\label{disrel}
     {\mathscr L}(x)=s_\lambda(\theta)\ge r_*(\theta-\theta_0).
    \eeq 

    \item There exists a number $L_*\ge 0$ such that for any $\lambda\in[0,1]$ and $\theta\ge \theta_0$,
    \beq\label{disrel2}
    {\mathscr L}(x)=s_\lambda(\theta)\le L_*+(\bar R+\mu_0)(\theta-\theta_0).
    \eeq

    \item  For any $\varep>0$, there exist numbers $\Theta_\varep\ge \Theta$ and $L_\varep\ge 0$ such that, for any $\lambda\in[0,1]$ and $\theta\ge \Theta_\varep$,
    \beq\label{disrel3}
    {\mathscr L}(x)=s_\lambda(\theta)\le L_\varep+(r_*+\mu_*+\varep)(\theta-\Theta_\varep).
    \eeq
\end{enumerate}
\end{lemma}
\begin{proof}
Let $\lambda$ be any given number in $[0,1]$.

(i) With $X_\lambda(\theta)=R_\lambda(\theta)(\cos(\theta),\sin(\theta))$, we have, for almost all $\theta\in[\theta_0,\infty)$, 
\beqs
X_\lambda'(\theta)=R_\lambda'(\theta)(\cos(\theta),\sin(\theta))+R_\lambda(\theta)(-\sin(\theta),\cos(\theta)),
\eeqs
thus,
\beq\label{Xpr}
|X_\lambda'(\theta)|=\sqrt{R_\lambda^2(\theta)+|R_\lambda'(\theta)|^2}
\le R_\lambda(\theta)+|R_\lambda'(\theta)|.
\eeq
Observe from the first identity in \eqref{Xpr} that $|X_\lambda'(\eta)|\ge R_\lambda(\eta)> r_*$, hence we have from \eqref{slform} the lower bound  \eqref{disrel}.

(ii) For $\theta\ge\theta_0$,  define
\beq\label{Lstar}
L_*(\theta)=\max\left\{ \int_{\theta_0}^{\theta} |r_1'(\eta)|\d\eta,\int_{\theta_0}^{\theta} |r_2'(\eta)|d\eta \right\},
\eeq 
which is finite thanks to \eqref{rloc}.
Considering  $\theta\ge \Theta$,  we estimate, by using the inequality in \eqref{Xpr}, the relations in \eqref{rrbound}, 
formula \eqref{Lstar} and the definition of $\mu_0$ in \eqref{muzero},
\begin{align*}
s_\lambda(\theta)
&\le \int_{\theta_0}^\theta R_\lambda(\eta)+|R_\lambda'(\eta)| \d \eta 
\le \bar R(\theta-\theta_0)+\int_{\theta_0}^\Theta |R_\lambda'(\eta)| \d \eta + \int_{\Theta}^\theta |R_\lambda'(\eta)| \d \eta\\
&\le \bar R(\theta-\theta_0)+L_*(\Theta)+\mu_0(\theta-\Theta)
= (\bar R+\mu_0)(\theta-\theta_0)+L_*(\Theta).
\end{align*}
Hence, we obtain \eqref{disrel2} with $L_*=L_*(\Theta)$.
Considering $\theta\in[\theta_0,\Theta)$ now, we have
\beqs
s_\lambda(\theta)
\le \bar R(\theta-\theta_0)+\int_{\theta_0}^\theta |R_\lambda'(\eta)| \d \eta
\le \bar R(\theta-\theta_0)+\int_{\theta_0}^\Theta |R_\lambda'(\eta)| \d \eta
\le \bar R(\theta-\theta_0) +L_*,
\eeqs
which yields \eqref{disrel2} again. Therefore, \eqref{disrel2} holds true for all $\theta\ge \theta_0$.

(iii) Let $\varep>0$. By  \eqref{bassum1} and definition \eqref{mustar}, there exists a number $\Theta_\varep\ge \theta_0$ such that, for all most all $\theta\ge \Theta_\varep$ and $i=1,2$,
\beq\label{rre}
r_i(\theta)\le r_*+\varep/2\text{ and }  |r_i'(\theta)|\le \mu_*+\varep/2.
\eeq 
This implies, for almost all $\eta\ge \Theta_\varep$, that
$R_\lambda(\eta)\le r_*+\varep/2$  and $|R_\lambda'(\eta)|\le \mu_*+\varep/2$, hence, 
\beq\label{rrsmall}
|X_\lambda'(\eta)|\le R_\lambda(\eta)+|R_\lambda'(\eta)|
\le r_*+ \mu_*+\varep.
\eeq
For $\theta\ge \Theta_\varep$, by rewriting $s_\lambda(\theta)$ and utilizing inequality \eqref{rrsmall}, we have
\beq\label{sl0}
s_\lambda(\theta)=s_\lambda(\Theta_\varep)+\int_{\Theta_\varep}^{\theta}|X_\lambda'(\eta)|\d\eta
\le s_\lambda(\Theta_\varep) +(r_*+\mu_*+\varep)(\theta-\Theta_\varep).
\eeq
Note from  \eqref{rrbound} and \eqref{Xpr}  that
\beq\label{sl1}
s_\lambda(\Theta_\varep)
\le \int_{\theta_0}^{\Theta_\varep} (\bar R + (1-\lambda)|r_1'(\eta)|+\lambda|r_2'(\eta)|)\d\eta\le L_\varep\eqdef \bar R(\Theta_\varep-\theta_0)
+L_*(\Theta_\varep).
\eeq
Thus, inequality \eqref{disrel3} follows from \eqref{sl0} and \eqref{sl1}.
\end{proof}

For each $\lambda\in[0,1]$, observe that $s_\lambda(\theta)$ is continuous in $\theta$, $s_\lambda(\theta_0)=0$ and, thanks to inequality \eqref{disrel}, $\lim_{\theta\to\infty}s_\lambda(\theta)=\infty$. Hence, 
\beqs 
s_\lambda([\theta_0,\infty))=[0,\infty), \text{ which implies }{\mathscr L}(X_\lambda([\theta_0,\infty))=[0,\infty). 
\eeqs 
Consequently, 
${\mathscr L}(\widetilde U)=[0,\infty)$.
For any subset $J$ of $[0,\infty)$ and number $\ell\ge 0$, denote
\beqs
\widehat U_J={\mathscr L}^{-1}(J)\text{ and  }
\widehat{\mathcal S}_\ell = \widehat U_{\{\ell\}} ={\mathscr L}^{-1}(\{\ell\})\ne \emptyset.
\eeqs
In particular, thanks to \eqref{slform} and \eqref{disrel}, $\widehat{\mathcal S}_0=\mathcal S_{\theta_0}$.
We now have two ways to view $\widetilde U$
\beqs
\widetilde U=\bigcup_{\theta\ge \theta_0}\mathcal S_\theta\text{ and }
\widetilde U=\bigcup_{\ell\ge 0}\widehat{\mathcal S}_\ell. 
\eeqs
These different decompositions  of $\widetilde U$ will be reflected in Theorem \ref{thm1} and Corollary \ref{cor1} below.

\begin{example}
We present some examples from simple to more complicated ones. Fix the numbers $r_*>0$, $\theta_0\in\R$ and $s>0$.
\begin{enumerate}[label=\rnum]
    \item For $ \theta\ge \theta_0$ and $i=1,2$,
\beqs 
r_i(\theta)=r_*\left\{1+\frac{1}{(\theta+z_i)^s}\right\},
 \text{ with }
z_2>-\theta_0,\ 
0\le z_2<z_1<z_2+2\pi.
\eeqs 

\item More generally, for $i=1,2$,
\beqs 
r_i(\theta)=r_*\left\{1+\frac{1}{(\theta+z_i+\varep_i\cos(\omega_i\theta))^s}\right\},\text{ or }
r_i(\theta)=r_*\left\{1+\frac{1}{(\theta+z_i+\varep_i\sin(\omega_i\theta))^s}\right\},
\eeqs 
with 
\beq\label{egcond}
z_1\ge 0,\ z_2\ge 0,\  
\theta_0+z_i-|\varep_i|>0,\ 
z_1-|\varep_1|>z_2+|\varep_2|,\ 
z_1+|\varep_1|<z_2-|\varep_2|+2\pi.
\eeq
The last two conditions can be simply rewritten as
\beqs 
|\varep_1|+|\varep_2|<z_1-z_2 < 2\pi-(|\varep_1|+|\varep_2|).
\eeqs
To verify \eqref{rnin}, observe that
\beqs 
r_*\left\{1+\frac{1}{(\theta+z_i+|\varep_i|)^s}\right\}\le r_i(\theta)\le r_*\left\{1+\frac{1}{(\theta+z_i-|\varep_i|)^s}\right\}.
\eeqs

\item Consider $\theta_0\ge 0$. For $i=1,2$,
\beqs 
r_i(\theta)=r_*\left\{1+\frac{1}{(\theta+z_i+\varep_i\cos(\omega_i\theta^{s+2}))^s}\right\},\text{ or }
r_i(\theta)=r_*\left\{1+\frac{1}{(\theta+z_i+\varep_i\sin(\omega_i\theta^{s+2}))^s}\right\},
\eeqs 
with the same conditions as in \eqref{egcond}.
\end{enumerate}

One can verify that  all examples {\rm (a)},  {\rm (b)} and {\rm (c)} above satisfy \eqref{rnin}, \eqref{bassum1} and Assumption \ref{bassum2}.
Moreover, in {\rm (a)} and {\rm (b)}, we have $r'_1(\theta)$ and $r'_2(\theta)$ go to $0$ as $\theta\to\infty$. Meanwhile,   we only have $r_1',r_2'\in L^\infty(\theta_0,\infty)$  in   {\rm (c)}.
\end{example}

\section{Growth lemmas}\label{glsec}

Consider the general $n$-dimensional space with $n\ge 1$. 
Denote the Euclidean norm of a vector $x\in\R^n$ by $|x|$.
Let $\mathcal M^{n\times n}_{{\rm sym}}$ denote the set of $n\times n$ symmetric matrices of real numbers.

Let $\Omega$ be an open subset of $\R^n$. 
Let $A=(a_{ij})_{i,j=1,\ldots,n}:\Omega\to \mathcal M^{n\times n}_{{\rm sym}}$ and $b:\Omega\to \R^n$  be given functions.
Define the linear differential operator $L$ by 
 \beq \label{Ldef} 
 Lu= -\sum_{i,j=1}^n a_{i,j}\frac{\partial^2 u}{\partial x_i \partial x_j}+b\cdot\nabla u
 \text{ for any function $u\in C^2(\Omega)$.}
 \eeq 
 
\subsection{A general Growth Lemma}

We recall a barrier function from \cite[Chapter 1, Lemma 2.1]{LandisBook}, see also  \cite[Chapter 1, Lemma 2.1]{Landis1968}.

\begin{lemma}\label{subsol}
Assume
\beq\label{welip}
\xi^{\rm T} A(x)\xi >0 \text{ for all $x\in \Omega$ and all }\xi\in \R^n\setminus\{0\},
\eeq
and there exist a point $x_*\not\in \Omega$ and a number $e_0\in \R$ such that 
\beq\label{b-cond}
\frac{{\rm Tr}A(x)-b(x)\cdot (x-x_*)} {\frac{(x-x_*)^{\rm T}}{|x-x_*|}A(x)\frac{(x-x_*)}{|x-x_*|}} \le e_0 \text{ for all }x\in \Omega.
\eeq   
Then the function $u(x)=|x-x_*|^{-s}$, for any positive number $s\ge e_0-2$, satisfies $u\in C^2(\Omega)$ and $Lu\le 0$ in $\Omega$.
\end{lemma}

The proof of Lemma \ref{subsol} is elementary and is given in Appendix \ref{apx}, with slightly general computations, for the sake of completeness.
The following two basic assumptions will often be considered throughout the paper.

\begin{assumption}[Ellipticity]\label{firstA}
    There are constants $c_0>0$ and $M_1>0$ such that
    \beqs
\xi^{\rm T} A(x)\xi \ge c_0|\xi|^2\text{ for all $x\in \Omega$ and all $\xi\in \R^n$,}
\eeqs  
 and 
\beqs
{\rm Tr}(A(x))\le M_1 \text{ for all }x\in \Omega.
\eeqs
\end{assumption}

\begin{assumption}[Bounded drift]\label{condB}
There is a constant $M_2\ge 0$ so that
\beqs\label{bM2}
 |b(x)|\le M_2 \text{ for all }x\in \Omega.
\eeqs
\end{assumption}

In the case $\Omega$ is bounded, we denote, for any point $y\not\in \bar\Omega$,
\beq\label{Rstar}
R_*(y,\Omega)=\max\{|x-y|:x\in\bar \Omega\}.
\eeq

For a real-valued function $f$, we recall that the positive and negative parts of $f$ are  $f^+=\max\{0,f\}$ and $f^-=\max\{0,-f\}=(-f)^+$. Then one has    
\begin{align*}
f=f^+-f^-,\  |f|=f^+ + f^-=\max\{f^+,f^-\},\\ 
0\le f^+\le |f|,\  0\le f^-\le |f|,
\text{ and } 
-f^-\le f\le f^+.
\end{align*} 
If a function $f:S\to \R$ is bounded from above, then
$\displaystyle \max\Big\{0,\sup_S f\Big\}=\sup_S (f^+)$.

\begin{lemma}[A general Growth Lemma]\label{GGL0}
Suppose $\Omega$ is bounded, 
\beqs
\partial \Omega=\widetilde\gamma\cup {\widetilde\Gamma} \text{ with }{\widetilde\Gamma}\ne \emptyset,
\eeqs 
and $S$ is a non-empty subset of $\bar \Omega$.
Assume there is a point $x_*\not\in \bar \Omega$ such that
\beq\label{Ddcond}
D_S\eqdef \sup_{x\in S}|x-x_*|< d_{\widetilde\Gamma}\eqdef \inf_{x\in {\widetilde\Gamma}} |x-x_*|.
\eeq
Let $d_*$ be a number such that
\beq\label{dminx}
0<d_*\le \min_{x\in\bar \Omega}|x-x_*|.
\eeq

Under Assumptions \ref{firstA} and \ref{condB}, let $u\in C^2(\Omega)\cap C(\bar \Omega)$ satisfy $Lu\le  0 $ in $\Omega$  and $u\le 0$ on $\widetilde\gamma$.
Then one has
\beq\label{growthu}
\sup_{x\in S} u^+(x)\le \eta_*\sup_{x\in {\widetilde\Gamma}} u^+(x),
\eeq 
where the number $\eta_*$ belongs to the interval $(0,1)$ and is defined by
\beq\label{etao}
\eta_*=1-\left\{\left (\frac{d_*}{D_S}\right)^s-\left(\frac{d_*}{d_{\widetilde\Gamma}}\right)^s \right\},
\eeq
with $s$ being any number such that 
\beq\label{ns2}
s>0\text{ and }s\ge \frac{M_1+M_2 R_*(x_*,\Omega)}{c_0}-2.
\eeq
\end{lemma}
\begin{proof}
Let $M=\max_{x\in \bar \Omega} u^+(x)\ge 0$.
Define, for $t>0$ and $x\in \bar \Omega$,
\beqs
\varphi(t)=(d_*/t)^s\text{ and }  V(x)=\varphi(|x-x_*|)= (d_*/|x-x_*|)^s.
\eeqs
Since $|x-x_*|\ge d_*$ for $x\in\bar \Omega$, one has $0\le V(x)\le 1$ in $\bar \Omega$.
With $D_S<d_{\widetilde\Gamma}$ in assumption \eqref{Ddcond}, we have
\beq\label{phph}
0<\varphi(d_{\widetilde\Gamma})<\varphi(D_S)\le 1.
\eeq

We apply Lemma \ref{subsol}, noting that we can take
\beqs
e_0=\frac{M_1+M_2 R_*(x_*,\Omega)}{c_0}
\eeqs 
in condition \eqref{b-cond}, and, hence, any number $s$ as in \eqref{ns2}. 
It results in $LV\le 0$ in $\Omega$.

Define, for $x\in \bar \Omega$, the function
\beq\label{Wdef}
W(x)=M(1-V(x)+\varphi(d_{\widetilde\Gamma}))\ge 0.
\eeq
Then $LW=-M LV\ge 0$ in $\Omega$.   
For $x\in \widetilde\gamma$, it follows from \eqref{Wdef} that 
\beq\label{Wub1}
W(x)\ge 0\ge u(x).
\eeq
For $x\in {\widetilde\Gamma}$, one has $|x-x_*|\ge d_{\widetilde\Gamma}$ and, hence,
\beq\label{Wub2}
W(x)=M(1-V(x)+\varphi(d_{\widetilde\Gamma}))\ge M(1-\varphi(d_{\widetilde\Gamma})+\varphi(d_{\widetilde\Gamma}))=M\ge u(x).
\eeq
With \eqref{Wub1} and \eqref{Wub2}, we have $W(x)\ge u(x)$ on the boundary $\partial \Omega$. This and the fact $Lu\le 0\le LW$ in $\Omega$ imply, by the Comparison Principle, that
\beq\label{uW}
u(x)\le W(x) \text{ in } \bar \Omega.
\eeq
Particularly, for $x\in S$, by \eqref{uW} and the fact $|x-x_*|\le D_S$, one has
\beq\label{ue1}
u(x)\le W(x)=M(1-V(x)+\varphi(d_{\widetilde\Gamma})\le M(1-\varphi(D_S)+\varphi(d_{\widetilde\Gamma}))=\eta_* M.
\eeq
Note that the number $\eta_*$ given in \eqref{etao} belongs to the interval $(0,1)$ thanks to \eqref{phph}. By the Maximum Principle and the fact $u\le 0$ on $\widetilde \gamma$, one has
\beq\label{MuG}
M= \max\{0,\max_{x\in \bar\Omega} u(x)\}= \max\{0,\max_{x\in \partial \Omega} u(x)\}= \max_{x\in {\widetilde\Gamma}} u^+(x).
\eeq
Therefore, inequalities \eqref{ue1} and \eqref{MuG} imply \eqref{growthu}.
\end{proof}

\subsection{Planar domains inside an annulus}
We return to the planar domain $U$ in section \ref{domsec}.
Given numbers $\bar\theta\ge \theta_0$, $d_0>0$ and $\theta_*\in (0,\pi/2)$. Assume that
\beq\label{inside}
r_2(\theta)\le r_*+d_0\text{ for all }\theta\in[\bar\theta,\bar\theta+2\theta_*], 
\eeq
and
\beq\label{coscond}
\cos\theta_*<\frac{r_*}{r_*+d_0}.
\eeq

 The following version of Growth Lemma is fundamental for studying the elliptic problems in our particular domain.
 
\begin{lemma}\label{GLthick}
Let
\beq\label{newOm}
\Omega=U_{(\bar \theta,\bar \theta +2\theta_*)} 
\text{ and  }
\widetilde\gamma=\gamma\cap \widetilde U_{[\bar \theta,\bar \theta +2\theta_*]}.
\eeq 
Under Assumptions \ref{firstA} and \ref{condB} for this domain $\Omega$,
suppose $Lu\le 0$ in $\Omega$ and $u\le 0$ on $\widetilde\gamma$.
Then there exists a number $\eta_*\in(0,1)$ independent of $u$ such that 
\beq\label{gineq}
\max_{x\in \mathcal S_{\bar\theta+\theta_*}} u^+(x)\le \eta_*\max \left\{ \max_{x\in \mathcal S_{\bar\theta}} u^+(x),\max_{x\in \mathcal S_{\bar\theta+2\theta_*}} u^+(x)\right\}.
\eeq 
More precisely, let $d_*$ and $s$ be any numbers such that
\beq\label{decond}
d_*> \max\left\{ (r_*+d_0)\left(\frac1{\cos\theta_*}-1\right), 
\frac{d_0^2}{2[(r_*+d_0)(1-\cos\theta_*)-d_0]} \right\},
\eeq
\beq\label{ns3}
s>0\text{ and }s\ge \frac{M_1+M_2(r_*+d_0+d_*)}{c_0}-2.
\eeq
Then one can take 
\beq \label{estar}
\eta_*=1-\left\{\left( \frac{d_*}{d_*+d_0}\right)^s- \left( \frac{d_*}{\widehat d_*}\right)^s\right\},
\eeq 
where 
\beq\label{hatd}
\widehat d_*=\sqrt{d_*^2+2d_*(r_*+d_0)(1-\cos\theta_*)} >d_*+d_0.
\eeq
\end{lemma}
\begin{proof}
We ignore $d_*$ in \eqref{decond} momentarily.
For $i=0,1,2$, let $\Gamma_i=\mathcal S_{\bar\theta+i\theta_*}$.
We apply  Lemma \ref{GGL0} to $\Omega$ and $\widetilde\gamma$ in \eqref{newOm}, together with  
$\widetilde\Gamma=\Gamma_0\cup \Gamma_2$ and $S=\Gamma_1$. 
The point $x_*$ and number $d_*$ in Lemma \ref{GGL0}  will be selected below, see Fig. \ref{PicThick}. 
    \begin{figure}[ht]
    \centering
    \includegraphics[scale=.4]{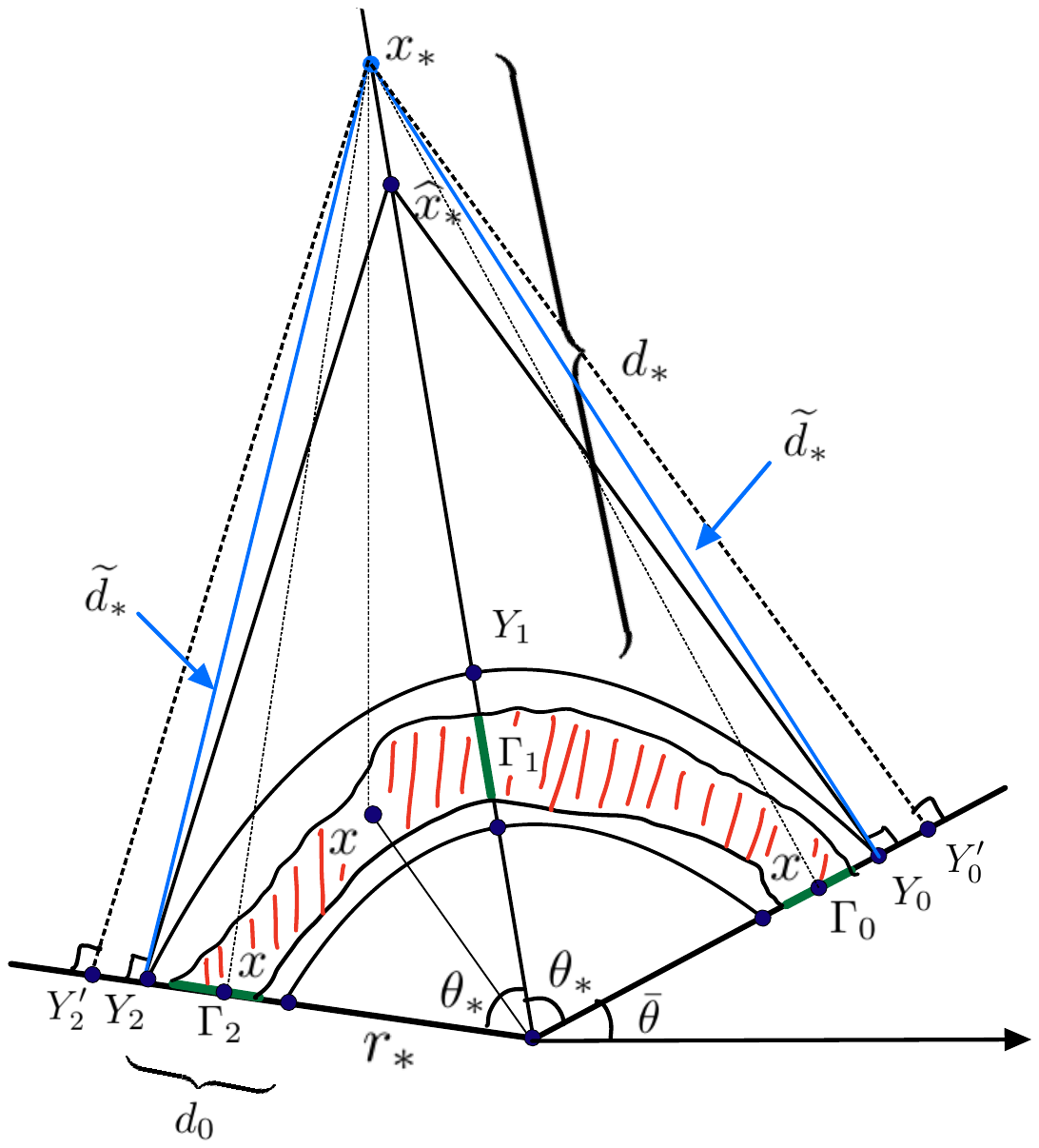}
    \caption{Construction for the Growth Lemma.}
    \label{PicThick}
 \end{figure}  
 
Let $\widetilde R=r_*+d_0$.
For $i=0,1,2$, denote
$Y_i=X(\widetilde R,\bar\theta+i\theta_*)$. 
We choose $x_*$ on the ray $X((\widetilde R,\infty)\times\{\bar \theta+\theta_*\})$ and set 
\beq \label{dstar}
d_*=|x_*-Y_1|=|x_*|-\widetilde R, \text{ that is, }|x_*|=d_*+r_*+d_0.
\eeq 
Let $\widehat x_*$ be the intersection of the tangent line to the circle $\{x\in\R^2:|x|=\widetilde R\}$ at the point $Y_0$ (or, $Y_2$) and the ray $\{X(r,\bar\theta+\theta_*):r>0\}$.
We take $x_*$ to be further than this $\widehat x_*$, that is, 
\beq\label{dxx}
|x_*|\ge |\widehat  x_*|=\frac{r_*+d_0}{\cos\theta_*}.
\eeq
By using  the last relation in \eqref{dstar}, we rewrite \eqref{dxx} as a condition on $d_*$ in the following way
\beq\label{decos}
d_*\ge \frac{r_*+d_0}{\cos\theta_*}-r_*-d_0=(r_*+d_0)\left(\frac1{\cos\theta_*}-1\right).
\eeq
In fact, the condition \eqref{decos} is met thanks to our assumption \eqref{decond}.

\medskip\noindent\emph{Verification of \eqref{dminx}.}
Thanks to \eqref{inside} and \eqref{dstar}, we have, for any $x=X(r,\theta)\in\bar\Omega$, 
\beqs
d_*=|x_*|-\widetilde R\le |x_*|-r_2(\theta)\le |x_*|-|x|\le |x-x_*|,
\eeqs 
which verifies \eqref{dminx}. 

\medskip\noindent\emph{Verification of \eqref{Ddcond}.} On the one hand, observe that
$\Gamma_1$ lies within the line segment $[0,x_*]$, hence
\beq \label{Ddcal1}
D_S=\max_{x\in\Gamma_1}|x-x_*|=|x_*|-r_1(\bar\theta+\theta_*)\le |x_*|-r_*= d_*+d_0.
\eeq 
On the other hand, set $\widetilde d_* =|x_*-Y_0|=|x_*-Y_2|$.
For $i=0,2$, letting $Y_i'$ be the orthogonal projection of $x_*$ onto the ray $\{X(r,\bar\theta+i\theta_*):r>0\}$. If $x\in\Gamma_i$, 
then $|x-Y_i'|\ge |Y_i-Y_i'|$, hence, using the right triangles $\Delta x_*Y_i'x$ and $\Delta x_*Y_i'Y_i$,  we have
$|x-x_*|\ge |Y_i-x_*|=\widetilde d_*$.
Therefore,
\beq \label{Ddcal2}
d_{\widetilde\Gamma}=\inf_{x\in \Gamma_0\cup\Gamma_2}|x-x_*|\ge \widetilde d_*.
\eeq 
We calculate $\widetilde d_*$ now. Recalling $|Y_0|=|Y_2|=\widetilde R$, we have 
\begin{align}
\notag
\widetilde d_*^2
&=|x_*|^2 +\widetilde R^2-2|x_*| \widetilde R\cos\theta_*=(|x_*|-\widetilde R)^2+2|x_*| \widetilde R(1-\cos\theta_*)\notag \\
&=d_*^2+2|x_*| \widetilde R(1-\cos\theta_*)>d_*^2+2d_* \widetilde R(1-\cos\theta_*)=\widehat d_*^2.    \label{dhati}
\end{align}
(The last inequality used the fact $|x_*|>d_*$.)
Combining \eqref{Ddcal2} with \eqref{dhati} yields 
\beq\label{Ddcal3}
d_{\widetilde\Gamma}\ge \widetilde d_* > \widehat d_*.
\eeq
Thanks to \eqref{Ddcal1} and \eqref{Ddcal3}, a particular sufficient condition for \eqref{Ddcond} is  
\beq \label{strictdd}
(d_*+d_0)^2< {\widehat d_*}^2,
\eeq 
which is   equivalent to
\beq \label{de}
    2d_* d_0+d_0^2< 2d_* \widetilde R(1-\cos\theta_*)=2d_* (r_*+d_0)(1-\cos\theta_*).
\eeq
Because of the assumption \eqref{coscond}, we have
\beqs 
(r_*+d_0)(1-\cos\theta_*)-d_0
=(r_*+d_0)\left[\frac{r_*}{r_*+d_0}-\cos\theta_*\right]>0.
\eeqs 
Thus, the sufficient condition \eqref{de} is equivalent to
\beqs
d_*> \frac{d_0^2}{2[(r_*+d_0)(1-\cos\theta_*)-d_0]},
\eeqs
which is satisfied thanks to \eqref{decond}.
Therefore, condition \eqref{Ddcond} is met.

\medskip\noindent\emph{Estimation of $R_*(x_*,\Omega)$.}
Consider $x=X(r,\theta)\in \bar\Omega$ with $\theta\in[\bar\theta,\bar\theta+2\theta_*]$. If $x$ is on the line segment $[0,x_*]$, then clearly
\beq \label{Rxx}
|x-x_*|\le |x_*|.
\eeq 
Otherwise, $\measuredangle{0xx_*}\ge \measuredangle{0Y_0 x_*}\ge \measuredangle{0Y_0\widehat x_*}=\pi/2$ in the case $\theta<\bar\theta+\theta_*$,
or $\measuredangle{0xx_*}\ge \measuredangle{0Y_2 x_*}\ge \measuredangle{0Y_2 \widehat x_*}=\pi/2$ in the case $\theta>\bar\theta+\theta_*$;
hence, from the triangle $\triangle 0xx_*$, we have \eqref{Rxx} again.
Therefore, recalling  \eqref{Rstar} and \eqref{dstar}, one has
\beq \label{Rxx2}
R_*(x_*,\Omega)=\max_{x\in\bar\Omega}|x-x_*|\le |x_*|=r_*+d_0+d_*.
\eeq 

\medskip Now, we can apply Lemma \ref{GGL0} and  obtain from \eqref{growthu} that 
\beq\label{gineq2}
\max_{x\in \mathcal S_{\bar\theta+\theta_*}} u^+(x)\le \eta\max \left\{ \max_{x\in \mathcal S_{\bar\theta}} u^+(x),\max_{x\in \mathcal S_{\bar\theta+2\theta_*}} u^+(x)\right\},
\eeq
where   the number $\eta$ can take the formula, thanks to \eqref{etao}, \eqref{ns2} and \eqref{Rxx2},
\beqs
\eta=1-\left\{\left (\frac{d_*}{D_S}\right)^s-\left(\frac{d_*}{d_{\widetilde\Gamma}}\right)^s \right\} \text{ with $s$ being in \eqref{ns3}.}
\eeqs
Thanks to the facts $D_S\le d_*+d_0$ in \eqref{Ddcal1} and $d_{\widetilde\Gamma}>\widehat d_*$ in \eqref{Ddcal3}, and $\eta_*$ being defined in \eqref{estar}, we have $\eta\le \eta_*$. Hence, inequality \eqref{gineq} follows from \eqref{gineq2}. Note that $\eta_*\in(0,1)$ thanks to \eqref{strictdd}.
\end{proof}

\section{The homogeneous problems}\label{estsec}

Hereafter, the operator $L$ is defined in \eqref{Ldef} under  Assumption \ref{firstA}  for $\Omega=U$.
In this section, we estimate the sub-solutions and solutions of the homogeneous problems.
Let $u\in C^2(U)\cap C(\widetilde U)$. Note that $u$ \emph{is not} assumed to be continuous on $\bar U$. In fact, it is not specified on the circle $\mathcal C_*$ which is a part of $\partial U$.

Let $\bar R$ be as in \eqref{rrbound} and $d_0=\bar R-r_*$. Assume $\theta_*\in(0,\pi/2)$ satisfies \eqref{coscond}. 
 For $i\ge 0$, set
 \beqs
\bar m_i = \max\big\{0, \max_{x\in  \mathcal S_{\theta_0+i\theta_*}} u(x) \big\}
=\max_{x\in \mathcal S_{\theta_0+i\theta_*}} u^+(x).
\eeqs

The following discrete dichotomy works for both bounded and unbounded drifts, see \cite[Theorem 6.1]{Landis1963} for elliptic equations, and also \cite[Lemma V.6]{HIK2} for a parabolic version in application. 

\begin{lemma}[Discrete dichotomy for solutions]  \label{cyls}
Assume 
\beq \label{blocal}
\text{$b(U_I)$ is bounded for any compact subsets 
$I$ in $[\theta_0,\infty)$.}
\eeq 
Let $d_*$ be as in \eqref{decond} and $\widehat d_*$ be defined by \eqref{hatd}.
 For $i\ge 1$, let
 \beqs
 B_i=\sup\{|b(x)|:x\in U_{\left(\theta_0+(i-1)\theta_*,\theta_0+(i+1)\theta_*\right)}\},
 \eeqs
 and fix any number $s_i$ such that
 \beq\label{si}
 s_i\ge \frac{M_1+B_i(r_*+d_0+d_*)}{c_0},
 \eeq
and define
 \beq\label{eti}
 \eta_i=1-\left\{\left( \frac{d_*}{d_*+d_0}\right)^{s_i}- \left( \frac{d_*}{\widehat d_*}\right)^{s_i}\right\}.
 \eeq

If $Lu\le 0$ in $U$ and $u\le 0$ on $\gamma$, then there are only the following two possibilities.
\begin{enumerate}[label=\tnum]
 \item \label{slemi} One has 
 \beq \label{ss0}
 \bar m_i\le \eta_i \bar m_{i-1}\text{ for all $i\ge 1$.} 
 \eeq
  Consequently, 
 \beq\label{ss1}
 \bar m_{i}\le \eta_{1}\eta_{2}\ldots \eta_{i}\bar m_{0}\text{ for all $i\ge 1$,}
 \eeq 
 \beq\label{ss2}
 \max_{x\in \mathcal S_\theta} u^+(x)
 \le  \max_{x\in \mathcal S_{\theta_0}} u^+(x)
 \text{ for all $\theta\ge \theta_0$.}
 \eeq
 
  \item \label{slemii} There is $i_0\ge 1$ such that 
 \beq \label{ss3}
 \bar m_{i_0}>0 \text{ and }\bar m_{i}\ge \frac{\bar m_{i-1}}{\eta_{i-1}} 
  \text{ for all $i\ge i_0+1$.}
 \eeq 
 Consequently,
\beq \label{ss4}
 \bar m_{i} 
 \ge \frac{\bar m_{i_0}}{\eta_{i_0}\eta_{i_0+1}\ldots \eta_{i-2}\eta_{i-1}}>0
 \text{ for all $i\ge i_0+1$.}
 \eeq
\end{enumerate}
\end{lemma}
\begin{proof} 
For any $i\ge 1$, applying Lemma \ref{GLthick} to $\bar\theta=\theta_0+(i-1)\theta_*$, one has from \eqref{gineq} that
\beqs
\bar m_{i}\le \eta_i\max\{\bar m_{i-1},\bar m_{i+1}\}.
\eeqs
Applying Lemma \ref{maindiseq} to the sequences $a_i=\bar m_{i}$ for $i\ge 0$ and $\lambda_i=\eta_{i}$ for $i\ge 1$, we obtain the dichotomy: either statement  \ref{slemi} but with only \eqref{ss0} and \eqref{ss1} is true, or the whole statement \ref{slemii} is true.
It remains to prove \eqref{ss2}. In the case \ref{slemi}, for any $\theta\ge \theta_0$, there is $i\ge 0$ such that $\theta\in[\theta_0+i\theta_*,\theta_0+(i+1)\theta_*)$. By the Maximum Principle for the domain $U_{(\theta_0+i\theta_*,\theta_0+(i+1)\theta_*)}$ and \eqref{ss1}, we have
\beqs
 \max_{x\in \mathcal S_\theta} u^+(x)\le \max\{\bar m_i,\bar m_{i+1}\}\le \bar m_0
 \eeqs
 which proves \eqref{ss2}.
\end{proof}

\subsection{Bounded drifts}\label{bddrift}

First, we state our results on the decay/growth of the sub-solutions and solutions in terms of the angle $\theta$ when the drift is bounded.

\begin{theorem}\label{thm1}
Under Assumption \ref{condB} for $\Omega=U$, there are constants $C_*>0$ and $\nu>0$ such that the following statements hold true.
\begin{enumerate}[label=\tnum]
    \item\label{T1i} Suppose $Lu\le 0$ in $U$ and $u\le 0$ on $\gamma$. Then there are only two possibilities.
    \begin{enumerate}[label=\rnum]
        \item\label{T1ia} One has
        \beq\label{d0a}
        u^+(x)\le \max_{x\in \mathcal S_{\theta_0}} u^+(x) \text{ for all }x\in \widetilde U,
        \eeq 
        \beq\label{d1a}
        \max_{x\in \mathcal S_\theta} u^+(x)\le C_*\left(\max_{x\in \mathcal S_{\theta_0}} u^+(x)\right)e^{-\nu(\theta-\theta_0)}\text{ for any }\theta\ge\theta_0.
        \eeq
        \item\label{T1ib} There are $C_1>0$ and $\bar\theta_1$ both depending on $u$ such that 
        \beq\label{d1b}
          \max_{x\in \mathcal S_\theta} u^+(x)\ge C_1e^{\nu(\theta-\theta_0)} \text{ for all }\theta\ge \bar \theta_1.
        \eeq 
    \end{enumerate}

    \item\label{T1ii} Suppose $Lu=0$ in $U$ and $u=0$ on $\gamma$. Then there are only two possibilities.
    \begin{enumerate}[label=\rnum]
        \item\label{T1iia} One has
        \beq\label{dabs0}
        |u(x)|\le \max_{x\in \mathcal S_{\theta_0}} |u(x)| \text{ for all }x\in \widetilde U,
        \eeq 
        \beq\label{dabsa}
        \max_{x\in \mathcal S_\theta} |u(x)|\le C_*\left(\max_{x\in \mathcal S_{\theta_0}} |u(x)|\right)e^{-\nu(\theta-\theta_0)} \text{ for any }\theta\ge\theta_0.
        \eeq
        \item\label{T1iib} There are $C>0$ and $\bar\theta$ both depending on $u$ such that
        \beq\label{dabsb}
         \max_{x\in \mathcal S_\theta} |u(x)| \ge Ce^{\nu(\theta-\theta_0)} \text{ for all }\theta\ge \bar \theta.
        \eeq  
    \end{enumerate}
\end{enumerate}
\end{theorem}
\begin{proof}
We use the same notation as in Lemma \ref{cyls}.

\ref{T1i} Thanks to Assumption \ref{condB}, $B_i\le M_2$ for all $i\ge 1$. We apply Lemma \ref{cyls} to
\beqs
s_i=s\eqdef \frac{M_1+M_2(r_*+d_0+d_*)}{c_0} \text{ for all $i\ge 1$.}
\eeqs
Then all $\eta_i$ in \eqref{eti} are $\eta_i=\eta_*\eqdef 1-\left\{\left( \frac{d_*}{d_*+d_0}\right)^{s}- \left( \frac{d_*}{\widehat d_*}\right)^{s}\right\}\in(0,1)$ for all $i\ge 1$.
Applying Lemma \ref{cyls}, we have the only two possibilities below.

\medskip\noindent\textit{Case 1 which corresponds to Lemma \ref{cyls}\ref{slemi}}. We will prove that this gives the statement \ref{T1i}\ref{T1ia} of this theorem. First, we have from \eqref{ss1} that
\beq\label{Mke}
\bar m_i\le \eta_*^i \bar m_0\text{ for all }i\ge 0.
\eeq
For $\theta\ge\theta_0$, there is an integer $i\ge 0$ such that $\theta-\theta_0\in[i\theta_*,(i+1)\theta_*)$ which implies
\beq\label{thek}
(\theta-\theta_0)/\theta_*-1<i\le (\theta-\theta_0)/\theta_*.
\eeq
By the Maximum principle for $u$ in the domain $U_{(\theta_0+i\theta_*,\theta_0+(i+1)\theta_*)}$ and \eqref{Mke}, we have
\beqs
\max_{x\in \mathcal S_\theta} u^+(x)\le\max\{\bar m_i,\bar m_{i+1}\} \le \bar m_0 \eta_*^i .
\eeqs
On the one hand, this immediately yields $\max_{x\in \mathcal S_\theta} u^+(x) \le \bar m_0 $ which proves \eqref{d0a}.
On the other hand, using also  \eqref{thek}, we continue to estimate
\beqs
\max_{x\in \mathcal S_\theta} u^+(x)\le \bar m_0 \eta_*^{(\theta-\theta_0)/\theta_*-1}=C_*\bar m_0 e^{-\nu(\theta-\theta_0)},
\eeqs
where
$C_*=\eta_*^{-1}$  and $\nu=\theta_*^{-1}\ln(1/\eta_*)$.
Thus, we obtain \eqref{d1a}.

\medskip\noindent\textit{Case 2 which corresponds to Lemma \ref{cyls}\ref{slemii}}. We will prove that this gives the statement \ref{T1i}\ref{T1ib} of this theorem.
Let $i_0$ be as in \eqref{ss3}. 
For any sufficiently large $\theta$, there is an integer $i\ge i_0+1$ such that  $\theta\in[ \theta_0+i\theta_*,\theta_0+(i+1)\theta_*)$.
Applying the Maximum Principle to the domain $U_{(\theta_0+(i-1)\theta_*,\theta)}$, we have
\beq\label{miimax}
\bar m_i\le\max\left \{\bar m_{i-1}, \max_{x\in \mathcal S_\theta} u^+(x)\right\}.
\eeq
Since $\bar m_i>\bar m_{i-1}$, see \eqref{ss3}, it follows \eqref{miimax} that 
$\bar m_i\le\max_{x\in \mathcal S_\theta} u^+(x)$.
Together with \eqref{ss4}, this implies
\beq \label{growthok}
\max_{x\in \mathcal S_\theta} u^+(x)\ge \bar m_i
\ge \frac{\bar m_{i_0}}{\eta_*^{i-i_0}}.
\eeq 
Note that $\theta_0+(i+1)\theta_*\ge \theta$, hence  $i\ge (\theta-\theta_0)/\theta_*-1$. 
Using this fact in \eqref{growthok} yields 
\beqs 
\max_{x\in \mathcal S_\theta} u^+(x)\ge \eta_*^{i_0+1} \bar m_{i_0} (\eta_*^{-1})^{(\theta-\theta_0)/\theta_*}=C_1e^{\nu(\theta-\theta_0)},
\eeqs 
where $C_1=\eta_*^{i_0+1} \bar m_{i_0}>0$. Therefore, we obtain inequality \eqref{d1b}.

\ref{T1ii} Applying part \ref{T1i} to $u:=-u$, we have two cases for $u^-$.
\begin{enumerate}[label=\rnum]
        \item One has
               \beq\label{d20}
        u^-(x)\le \max_{x\in \mathcal S_{\theta_0}} u^-(x) \text{ for all }x\in U,
        \eeq 
        \beq\label{d2a}
        \max_{x\in \mathcal S_{\theta}} u^-(x)\le C_*\left(\max_{x\in \mathcal S_{\theta_0}} u^-(x)\right)e^{-\nu(\theta-\theta_0)} \text{ for all }\theta\ge \theta_0.
        \eeq
        \item  There are $C_2>0$ and $\bar\theta_2$ such that
        \beq\label{d2b}
           \max_{x\in \mathcal S_{\theta}} u^-(x)\ge C_2e^{\nu(\theta-\theta_0)} \text{ for all }\theta\ge \bar \theta_2.
        \eeq 
    \end{enumerate}
    Combining these with the two possibilities in part \ref{T1i}, all together they fall into three cases below.

    Case 1: \eqref{d0a}, \eqref{d1a}, \eqref{d20} and \eqref{d2a} hold. Then using $|u(x)|=\max\{u^+(x),u^-(x)\}$, we obtain 
    \eqref{dabs0} from \eqref{d0a} and \eqref{d20}, and obtain 
    \eqref{dabsa} from \eqref{d1a}, \eqref{d2a}.
    
    Case 2: \eqref{d1b} holds. Then using $|u(x)|\ge  u^+(x)$ and \eqref{d1b}, we obtain \eqref{dabsb}.

    Case 3. \eqref{d2b} holds. Then using $|u(x)|\ge u^-(x)$ and \eqref{d2b}, we obtain \eqref{dabsb}.
\end{proof}

If Assumption \ref{bassum2} is also true, then Theorem \ref{thm1} can be restated using the arc-distance ${\mathscr L}(x)$ .
Define 
\beq\label{kstar}
k_*=\frac{\bar R+\mu_0}{r_*}>1.
\eeq

\begin{corollary}\label{cor1}
Under Assumption  \ref{bassum2}  and Assumption \ref{condB} for $\Omega=U$,  there are positive numbers $C_*'>0$ and $\nu'>0$ such that the following statements hold true.
\begin{enumerate}[label=\tnum]
    \item\label{cori} Suppose $Lu\le 0$ in $U$ and $u\le 0$ on $\gamma$. Then there are only two possibilities.
    \begin{enumerate}[label=\rnum]
        \item\label{coria} One has \eqref{d0a} and 
        \beq\label{cs1}
        u^+(x)\le C_*'\left(\max_{x\in \mathcal S_{\theta_0}} u^+(x)\right)e^{-\nu' {\mathscr L}(x)} \text{ for all }x\in \widetilde U.
        \eeq
        \item \label{corib}
            For any $\varep>0$, there are numbers $C_1'>0$ depending on $u$  and $\bar\ell_1\ge 0$  depending on $\varep$ such that if $\ell\ge \bar\ell_1$, then
            \beq\label{cs2}
         \sup\left\{u^+(x):x\in\widehat U_{\left[(k_*+\varep)^{-1}\ell,(k_*+\varep)\ell\right]}\right\}\ge C_1'e^{\nu' \ell}.
        \eeq
    \end{enumerate}

    \item\label{corii} Suppose $Lu=0$ in $U$ and $u=0$ on $\gamma$. 
    Then there are only two possibilities.
    \begin{enumerate}[label=\rnum]
        \item\label{coriia} One has \eqref{dabs0} and 
        \beq\label{cs3}
        |u(x)|\le C_*'\left(\max_{x\in \mathcal S_{\theta_0}} |u(x)|\right)e^{-\nu' {\mathscr L}(x)} \text{ for all }x\in \widetilde U.
        \eeq
        \item \label{coriib}
        There is a number $C'>0$ depending on $u$  such that for any $\varep>0$,   there exists a number $\bar\ell\ge 0$  depending on $\varep$ so that if $\ell\ge \bar\ell$, then
               \beq\label{cs4}
        \sup\left\{|u(x)|:x\in\widehat U_{\left[(k_*+\varep)^{-1}\ell,(k_*+\varep)\ell\right]}\right\}\ge C'e^{\nu' \ell}.
        \eeq
    \end{enumerate}
\end{enumerate}
\end{corollary}
\begin{proof}
We apply Theorem \ref{thm1}.
Consider statements (a) of (i) and (ii) in Theorem \ref{thm1}. 
For $x=X(r,\theta)\in \widetilde U$, it follows from relation \eqref{disrel} and \eqref{disrel2} that
\beq\label{Lrx}
({\mathscr L}(x)-L_*)/(\bar R+\mu_0)  \le \theta-\theta_0\le {\mathscr L}(x)/r_*.
\eeq
Then
    \beq\label{etsx}
   e^{-\nu\left(\theta-\theta_0\right)}\le e^{-\nu ({\mathscr L}(x)-L_*)/(\bar R+\mu_0)}=e^{\nu L_*/(\bar R+\mu_0)} e^{-\nu'{\mathscr L}(x)},
   \text{ where }\nu'=\nu/(\bar R+\mu_0).
    \eeq
Using \eqref{etsx} in \eqref{d1a} and \eqref{dabsa}, we obtain \eqref{cs1} and \eqref{cs3}, respectively, with $C_*'=C_* e^{\nu L_*/(\bar R+\mu_0)}$.

Now, consider statements (b) of (i) and (ii) in Theorem \ref{thm1}. 
For $\ell\ge 0$, we pick a point $x\in\mathcal S_\theta$ for some $\theta\ge \theta_0$ such that ${\mathscr L}(x)=\ell$. Consider sufficiently large   $\ell$  and $\theta$.  Using \eqref{etsx}, we have
    \beq\label{etsx2}
   e^{\nu(\theta-\theta_0)}\ge e^{-\nu L_*/(\bar R+\mu_0)} e^{\nu'\ell}.
    \eeq

Let $y\in \mathcal S_\theta$. Combining \eqref{disrel2} for $x:=y$ with  \eqref{Lrx} gives
\beqs
{\mathscr L}(y)\le L_*+(\bar R+\mu_0)(\theta-\theta_0)
 \le L_*+(\bar R+\mu_0){\mathscr L}(x)/r_*=L_*+k_*\ell.
\eeqs
Combining \eqref{disrel} for $x:=y$ with  \eqref{Lrx} gives
\beqs
{\mathscr L}(y)\ge r_*(\theta-\theta_0)
\ge r_*({\mathscr L}(x)-L_*)/(\bar R+\mu_0)=k_*^{-1}\ell-k_*^{-1}L_*.
\eeqs
Therefore, $y\in \widehat U_{\left[k_*^{-1}\ell-k_*^{-1}L_*,k_*\ell+L_*\right]}$.
Consequently,
\beq\label{UU4}
\mathcal S_\theta\subset \widehat U_{\left[k_*^{-1}\ell-k_*^{-1}L_*,k_*\ell+L_*\right]}.
\eeq

For any $\varep>0$, we can take $\ell$ sufficiently large such that  
\beq\label{ke}
(k_*+\varep)^{-1}\ell\le k_*^{-1}\ell-k_*^{-1}L_*< k_*\ell+L_*\le (k_*+\varep)\ell, 
\eeq
Thanks to \eqref{disrel2}, $\theta$ is also sufficiently large. 
Then combining \eqref{UU4} with \eqref{ke} yields
$\mathcal S_\theta\subset  \widehat U_{\left[(k_*+\varep)^{-1}\ell,(k_*+\varep)\ell\right]}$.
Using this property for the left-hand sides of \eqref{d1b} and \eqref{dabsb}, and using \eqref{etsx2} for the right-hand sides of those same inequalities, we obtain 
\eqref{cs2} with $C_1'=C_1 e^{-\nu L_*/(\bar R+\mu_0)}$
and 
\eqref{cs4} with $C'=C e^{-\nu L_*/(\bar R+\mu_0)}$.
\end{proof}

One can obviously rewrite \eqref{cs1} as
\beq\label{Elim}
\sup_{x\in \widehat{\mathcal S}_\ell} u^+(x)\le C_*'\left(\max_{x\in \widehat{\mathcal S}_0} u^+(x)\right)e^{-\nu'\ell},\text{ hence, }
\lim_{\ell\to\infty} \left(\sup_{x\in \widehat{\mathcal S}_\ell} u^+(x) \right)=0.
\eeq
In the case of \eqref{cs3}, we obtain \eqref{Elim} with $|u(x)|$ replacing $u^+(x)$.

\begin{remark}\label{near1}
 Consider the case $\mu_*=0$. (Such a case occurs particularly when $r_1'(\theta)$ and $r_2'(\theta)$ exist for all sufficiently large $\theta$ and both go to zero as $\theta\to\infty$.)
 Let $\varep>0$ be sufficiently small such that $r_*+\varep/2<\bar R$, and the number
 $k_\varep\eqdef \frac{r_*+\varep}{r_*}+\varep$ is close to $1$.
 Consider large $\theta$ such that \eqref{rre} is satisfied.
In the proofs of \eqref{cs2} and \eqref{cs4} in Corollary \ref{cor1}, we can still use the same numbers $C_*$ and $\nu$, but use, instead of \eqref{disrel2}, the following convenient consequence of \eqref{disrel3}
   \beqs
    {\mathscr L}(x)=s_\lambda(\theta)\le L_\varep+(r_*+\mu_*+\varep)(\theta-\theta_0).
    \eeqs
Then we can replace $\mu_0$ with $\mu_*+\varep/2$, replace $\bar R$ with $r_*+\varep/2$, replace $L_*$ with $L_\varep$,   replace $k_*$ in \eqref{kstar} with $\frac{(r_*+\varep/2)+(\mu_*+\varep/2)}{r_*}=\frac{r_*+\varep}{r_*}$. Then estimates \eqref{cs2} and \eqref{cs4} hold when $\ell$ is large with $k_*+\varep$ being replaced with $k_\varep$ which is close to $1$. These are improved estimates for large $\ell$.
\end{remark}

To select which cases  to occur in Theorem \ref{thm1} and Corollary \ref{cor1}, we introduce the following class of functions with some restricted behavior as $\theta\to\infty$.

\begin{definition}\label{subex}
Let $K\in\{\widetilde U,U,\gamma_1,\gamma_2,\gamma\}$.    A function $f:K\to \R$ is said to be below exponential growth (in/on $K$) if, for any $\varep>0$,
\beqs
\lim_{\theta\to\infty} \left(e^{-\varep \theta} \sup_{x\in \mathcal S_\theta\cap K} |u(x)|\right) =0.
\eeqs 
\end{definition}

With Definition \ref{subex}, the following is an obvious consequence of Theorem \ref{thm1} and Corollary \ref{cor2}.

\begin{corollary}\label{cor2}
   Assume $u(x)$ is below exponential growth in $\widetilde U$. Then the statements  \ref{T1i}\ref{T1ia} and \ref{T1ii}\ref{T1iia} in Theorem \ref{thm1}, as well as \ref{cori}\ref{coria} and \ref{corii}\ref{coriia} in Corollary \ref{cor1} hold true.
\end{corollary}

\begin{remark}\label{similar}
    Above in this subsection \ref{bddrift}, we presented the results in terms of both $\theta$ and ${\mathscr L}(x)$. To avoid lengthening the paper, the remaining results below are only written for $\theta$ even though their counterparts for $\mathscr L(x)$ are also true.
\end{remark}

\subsection{Unbounded drifts}\label{ubdrift}
We derive estimates for the sub-solutions and solutions even when the drift is unbounded in $\widetilde U$.

\begin{assumption}\label{bmassum}
There are a number $\bar\Theta\ge \theta_0$ and an increasing function 
 $m:[\bar\Theta,\infty)\to [0,\infty)$ such that 
$b(x)$ is bounded in $U_{[\theta_0,\bar\Theta]}$ and,  for any $\theta\ge \bar\Theta$,  
\beqs
|b(x)|\le m(\theta) \text{ for all }x\in U\cap \mathcal S_\theta,
\eeqs
\beqs
\lim_{\theta\to\infty}m(\theta)=\infty,
\eeqs
and  there is a number $\kappa\in(0,1)$ so that
\beq\label{mt1}
\int_{\bar\Theta}^\infty \kappa^{m(\theta)}\d \theta=\infty.
\eeq
\end{assumption}

Thanks to the monotonicity of  the function $\theta\mapsto \kappa^{m(\theta)}$,
the integral in \eqref{mt1} is both improper Riemann integral and the Lebesgue integral.

\begin{theorem}\label{thm3}
Under Assumption \ref{bmassum}, there exist numbers $\Theta_0\ge \bar\Theta$ and $\bar\nu>0$ such that one has the following. If $Lu\le 0$ in $U$ and $u\le 0$ on $\gamma$, then there are only two possibilities.
\begin{enumerate}[label=\tnum]
    \item One has, for all $\theta\ge\theta_0$,
\beq\label{ubs1}
\max_{x\in \mathcal S_\theta} u^+(x)\le \left(\max_{x\in \mathcal S_{\theta_0}} u^+(x)\right) \min\left\{1,\exp\left( -\bar\nu\int_{\Theta_0}^\theta \kappa^{m(\zeta)} \d\zeta\right)\right\}.
\eeq
Consequently,
\beq\label{ubs2}
\lim_{\theta\to\infty} \left(\max_{x\in \mathcal S_{\theta}} u^+(x)\right)=0.
\eeq

\item There exists $C>0$ depending on $u$  such that, for sufficiently large $\theta$,
\beq\label{ubs3}
\max_{x\in \mathcal S_\theta} u^+(x) \ge C \exp\left( \bar\nu\int_{\bar\Theta}^\theta \kappa^{m(\zeta)} \d\zeta\right).
\eeq
Consequently,
\beq\label{ubs4}
\lim_{\theta\to\infty} \left(\max_{x\in \mathcal S_{\theta}} u^+(x)\right)=\infty.
\eeq
\end{enumerate}
\end{theorem}
\begin{proof}
We divide the proof into many steps.

\medskip\noindent\textit{Step 1.} Clearly, $b(x)$ satisfies condition \eqref{blocal}. We apply Lemma \ref{cyls} and refer to \ref{slemi} and \ref{slemii} there as Case Ia and Case Ib respectively.

\medskip\noindent\textit{Step 2.} We will improve the estimates for $u$ when $\theta$ is large. We will apply Lemma \ref{cyls} again with the set of parameters $(\theta_0,d_0,\theta_*,d_*)$ being replaced with $(\bar\theta_0,\bar d_0,\bar \theta_*,\bar d_*)$ below. We proceed as follows.

Firstly, let  $\bar\theta_*=\theta_*/2$ and pick $\varep_*\in(0,d_0)$ sufficiently small so that 
\beq \label{newecond}
\frac{1}{e^{\varep_*/c_0}}>\kappa\text{ and }
\cos\bar\theta_*<\frac{r_*}{r_*+\varep_*}.
\eeq

Secondly, set $\bar d_0=\varep_*$ and take an integer $N\ge 0$ sufficiently large such that $\bar\theta_0\eqdef \theta_0+N\theta_*$ satisfies
\beqs
\bar\theta_0\ge \bar\Theta\text{ and } r_2(\theta)\le r_*+\varep_*=r_*+\bar d_0 \text{ for all }\theta\ge \bar\theta_0.
\eeqs

Thirdly, choose $\bar d_*$ sufficiently large such that \eqref{decond} is satisfied   with 
\beq\label{replace}
\text{ $\bar d_*$ replacing $d_*$, $\bar d_0$ replacing $d_0$, and $\bar\theta_*$ replacing $\theta_*$, }
\eeq
and also
\beq\label{kkcond}
\bar\kappa\eqdef \frac{\bar d_*}{\bar d_*+\bar d_0}\in[\kappa,1)\text{ and }
\bar\kappa^\frac{r_*+\bar d_0+\bar d_*}{c_0}\ge \kappa.
\eeq
The last requirement is possible because, thanks to \eqref{newecond}, 
$$\lim_{\bar d_*\to\infty}\bar\kappa^\frac{r_*+\bar d_0+\bar d_*}{c_0}=e^{-\bar d_0/c_0}=e^{-\varep_*/c_0}>\kappa.$$

Fourthly, let $\widehat{\bar d_*}$ be defined as $\widehat d_*$ in \eqref{hatd} with the replacements in \eqref{replace}. For  $i\ge 1$, denote
\beq\label{newchoice}
\begin{aligned} 
\bar B_i&=m(\bar\theta_0+(i+1)\bar\theta_*),\quad 
\bar s_i= \frac{M_1+\bar B_i(r_*+\bar d_0+\bar d_*)}{c_0},\quad 
\bar z= \frac{\bar d_*+\bar d_0}{\widehat{\bar d_*}}\in(0,1),\\
\bar\eta_i&=1-\left[ \left(\frac{\bar d_*}{\bar d_*+\bar d_0}\right)^{\bar s_i}-\left(\frac{\bar d_*}{\widehat{\bar d_*}}\right)^{\bar s_i}\right]
    = 1- \bar\kappa^{\bar s_i}(1-\bar z^{\bar s_i}).
\end{aligned}
\eeq 

Fifthly, for $i\ge 0$, let $\displaystyle \bar M_i = \max_{x\in \mathcal S_{\bar\theta_0+i\bar\theta_*}} u^+(x)$.
We apply Lemma \ref{cyls} again this time with $\bar\theta_0$ replacing $\theta_0$, $r_*+\bar d_0$ replacing $\bar R$, and \eqref{replace}, and $\widehat{\bar d_*}$ replacing $\widehat d_*$ .
It results in the discrete dichotomy as in Lemma \ref{cyls} with $\bar M_i$ replacing $\bar m_i$, and $\bar s_i$ and $\bar \eta_i$ in \eqref{newchoice} replacing $s_i$ in \eqref{si} and $\eta_i$ in \eqref{eti}, respectively. We will refer to \ref{slemi} and \ref{slemii} in Lemma \ref{cyls} this time as Case IIa and Case IIb.

\medskip\noindent\textit{Step 3.} 
For $j\ge i\ge1$, we need to estimate 
$E_{i,j}\eqdef \prod_{k=i}^j \bar \eta_k$.
Since $\bar s_k$ is increasing in $k$, one has, for  $k\ge 1$, $1-\bar z^{\bar s_k}\ge \lambda_0\eqdef1-\bar z^{\bar s_1}$. 
Hence, $\bar\eta_k\le  1-\lambda_0\bar z^{\bar s_k}$, and we have 
\beq\label{Eks}
E_{i,j}\le \prod_{k=i}^j\left(1-\lambda_0\bar\kappa^{\bar s_k}\right)
\text{ which yields }
\ln E_{i,j}\le \sum_{k=i}^j \ln \left(1-\lambda_0\bar\kappa^{\bar s_k}\right)
\le -\lambda_0\sum_{k=i}^j \bar\kappa^{\bar s_k}.
\eeq
Denote
$\displaystyle S(\theta)=\frac{M_1+m(\theta)(r_*+\varep_*+\bar d_*)}{c_0}$.
Then
$$\bar s_k=S(\bar\theta_0+(k+1)\bar\theta_*)\le S(\zeta) 
\text{ for all }\zeta \in[\bar\theta_0+(k+1)\bar\theta_*,\bar\theta_0+(k+2)\bar\theta_*] .$$
This implies the last summand $\bar\kappa^{\bar s_k}$ in \eqref{Eks} is bounded from above by the average of $\bar\kappa^{S(\zeta)}$ over the last interval of $\zeta$. Therefore, we obtain
\beqs
\ln E_{i,j}
\le- \lambda_0 \sum_{k=i}^j\left(\frac{1}{\bar \theta_*}\int_{\bar \theta_0+(k+1)\bar \theta_*}^{\bar \theta_0+(k+2)\bar \theta_*} \bar\kappa^{S(\zeta)}\d\zeta\right)
=-  \frac{\lambda_0}{\bar \theta_*}  \int_{\bar \theta_0+(i+1)\bar \theta_*}^{\bar \theta_0+(j+2)\bar \theta_*} \bar\kappa^\frac{M_1}{c_0} \left(\bar\kappa^\frac{r_*+\bar d_0+\bar d_*}{c_0}\right)^{m(\zeta)}\d\zeta.
\eeqs
Thanks to condition \eqref{kkcond}, it follows that
\beqs
\ln E_{i,j}
\le -\bar\nu \int_{\bar\theta_0+(i+1)\bar\theta_*}^{\bar\theta_0+(j+2)\bar\theta_*} \kappa^{m(\zeta)}\d\zeta,
\text{ where }
\bar\nu=\frac{\lambda_0}{\bar \theta_*} \bar \kappa^{\frac{M_1}{c_0}}.
\eeqs
Hence, we have
\beq\label{Eibound} 
E_{i,j}
\le \exp\left(-\bar\nu \int_{\bar\theta_0+(i+1)\bar\theta_*}^{\bar\theta_0+(j+2)\bar\theta_*} \kappa^{m(\zeta)}\d\zeta\right),
\eeq
and, in particular when $i=1$ and $j\ge 1$,
\beq\label{E1bound}
E_{1,j}
\le \exp\left(-\bar\nu \int_{\Theta_0}^{\bar\theta_0+(j+2)\bar\theta_*} \kappa^{m(\zeta)}\d\zeta\right),
\text{ where } \Theta_0=\bar\theta_0+2\bar\theta_*=\theta_0+(N+1)\theta_*.
\eeq

\medskip\noindent\textit{Step 4.}
Consider the Case IIa. For all $i\ge 1$, we have
\beq\label{Munb1}
\bar M_i\le \bar \eta_i \bar M_{i-1}\text{ and }
\bar M_i\le \bar \eta_i \bar \eta_{i-1}\ldots \bar \eta_1 \bar M_0=E_{1,i}\bar M_0.
\eeq
Note that 
$\bar m_{N+i}=\bar M_{2i}\le \bar \eta_{2i}\bar \eta_{2i-1} \bar M_{2i-2}=\eta_{2i}\bar \eta_{2i-1}\bar m_{N+(i-1)},$
hence, Case Ia must also occur (not Case Ib).
In particular, 
\beq\label{Mmzero} 
\bar M_0=\bar m_N\le \bar m_0.
\eeq 

Consider $\theta\in[\theta_0,\Theta_0)$. Choose $i$ large such that $\theta\le \bar\theta_0+i\bar\theta_*$. By the Maximum Principle,
\beq\label{smallt}
\max_{x\in \mathcal S_\theta} u^+(x)\le \max\{ \bar m_0,\bar M_i\}\le \bar m_0.
\eeq

Consider $\theta\ge \Theta_0$ now. There is an integer $j\ge 1$ such that $\theta\in[\bar\theta_0+j\bar\theta_*,\bar\theta_0+(j+1)\bar\theta_*)$. By the Max Principle, \eqref{Munb1}, \eqref{Mmzero}  and  estimate \eqref{E1bound},
\begin{align}\notag
\max_{x\in \mathcal S_\theta} u^+(x)
&= \max\{\bar M_j,\bar M_{j+1}\}
= \bar M_j\le \bar m_0\exp\left(-\bar\nu \int_{\Theta_0}^{\bar\theta_0+(j+2)\bar\theta_*} \kappa^{m(\zeta)}\d\zeta\right)\\
&\le \bar m_0\exp\left(-\bar\nu \int_{\Theta_0}^{\theta} \kappa^{m(\zeta)}\d\zeta\right).
\label{larget}
\end{align}
The last inequality used the simple fact $\bar\theta_0+(j+2)\bar\theta_*>\theta$.

From \eqref{smallt} and \eqref{larget}, we obtain \eqref{ubs1} for all $\theta\ge\theta_0$. Then the  limit \eqref{ubs2} clearly follows from \eqref{mt1} and \eqref{ubs1}. 

\medskip\noindent\textit{Step 5.} 
Consider Case IIb. 
Then there exists an integer $i_0\ge 1$ such that $\bar M_{i_0}>0$ and, for $j\ge i_0+1$,
\beqs
\bar M_j\ge \frac{\bar M_{j-1}}{\bar\eta_{j-1}}
\text{ and }
\bar M_j\ge \frac{\bar M_{i_0}}{E_{i_0,j-1}}.
\eeqs
Combining the last inequality with estimate \eqref{Eibound} gives
\beqs
\bar M_j
\ge \bar M_{i_0}\exp\left(\bar\nu \int_{\bar\theta_0+(i_0+1)\bar\theta_*}^{\bar\theta_0+(j+1)\bar\theta_*} \kappa^{m(\zeta)}\d\zeta\right)
=C_2\exp\left(\bar\nu \int_{\bar\Theta}^{\bar\theta_0+(j+1)\bar\theta_*} \kappa^{m(\zeta)}\d\zeta\right).
\eeqs
where
\beqs
C_2=\bar M_{i_0}\exp\left(-\bar\nu \int_{\bar\Theta}^{\bar\theta_0+(i_0+1)\bar\theta_*} \kappa^{m(\zeta)}\d\zeta\right)>0.
\eeqs

Consider $\theta> \bar\theta_0+(i_0+2)\bar\theta_*$. There is  $j\ge i_0+2$ such that $\theta\in[\bar\theta_0+j\bar\theta_*,\bar\theta_0+(j+1)\bar\theta_*)$. 
Using the same arguments as in the proof of Theorem \ref{thm1}\ref{T1i}\ref{T1ib}, we have, see the first inequality in \eqref{growthok},
\begin{align*}
\max_{x\in \mathcal S_\theta} u^+(x)&\ge \bar M_j
\ge C_2\exp\left(\bar\nu \int_{\bar\Theta}^{\bar\theta_0+(j+1)\bar\theta_*} \kappa^{m(\zeta)}\d\zeta\right)
\ge C_2\exp\left(\bar\nu \int_{\bar\Theta}^{\theta} \kappa^{m(\zeta)}\d\zeta\right).
\end{align*}
Hence, we obtain \eqref{ubs3} and, together with \eqref{mt1}, it implies the limit \eqref{ubs4}.
\end{proof}

The estimates for solutions of the equation  are derived next. 

\begin{theorem}\label{thm4}
Under Assumption \ref{bmassum}, suppose  $Lu=0$ in $U$ and $u=0$ on $\gamma$.
Let $\Theta_0\ge \bar\Theta$ and $\bar\nu>0$ be the same numbers as in Theorem \ref{thm3}.
Then there are only two possibilities.
\begin{enumerate}[label=\tnum]
    \item\label{T4i} One has, for all $\theta\ge\theta_0$,
\beq\label{ubdd1}
\max_{x\in \mathcal S_\theta} |u(x)|\le \left(\max_{x\in \mathcal S_{\theta_0}}|u(x)|\right) \min\left\{1,\exp\left( -\bar\nu\int_{\Theta_0}^\theta \kappa^{m(\zeta)} \d\zeta\right)\right\},
\eeq
and, consequently,
\beqs
\lim_{\theta\to\infty} \left(\max_{x\in \mathcal S_\theta} |u(x)|\right)=0.
\eeqs
\item\label{T4ii} There exists $C>0$ depending on $u$  such that for sufficiently large $\theta$ 
\beqs
\max_{x\in \mathcal S_\theta} |u(x)|\ge C \exp\left( \bar\nu\int_{\bar\Theta}^\theta \kappa^{m(\zeta)} \d\zeta\right),
\eeqs
and, consequently,
\beqs
\lim_{\theta\to\infty} \left(\max_{x\in \mathcal S_\theta} |u(x)|\right)=\infty.
\eeqs
\end{enumerate}
\end{theorem}
\begin{proof}
We can prove this theorem based on Theorem \ref{thm3} in the same way that 
part \ref{T1ii} of Theorem \ref{thm1} was proved based on part \ref{T1i} of that same theorem. We omit the details.
\end{proof}

\section{Inhomogeneous problems}\label{inhomsec}

In this section, we consider only the bounded drifts. Therefore, we suppose, throughout the section,  Assumptions \ref{firstA} and \ref{condB} hold with $\Omega=U$. We assume additionally that  $u\in C^2(U)\cap C(\widetilde U)$  is \textit{bounded on $\gamma$} and is \textit{below exponential growth} in $\widetilde U$.
Define $M_\gamma=\sup_{x\in\gamma} u(x)$, $m_\gamma=\inf_{x\in\gamma} u(x)$,
and, for $\theta\ge \theta_0$,
\beqs
M_\gamma(\theta)=\max\{ u(X(r_i(\theta),\theta)): i=1,2\}, \quad
m_\gamma(\theta)=\min\{ u(X(r_i(\theta),\theta)):i=1,2\}.    
\eeqs 
Note that $M_\gamma-m_\gamma$ is the oscillation of $u$ on $\gamma$, and $\displaystyle\limsup_{\theta\to\infty}M_\gamma(\theta)-\liminf_{\theta\to\infty}m_\gamma(\theta)$ can be seen as the asymptotic  oscillation of $u$ on $\gamma$.

Below, $C_*$ and $\nu$ are the positive constants in Theorem \ref{thm1}.
Different scenarios for the inhomogeneity are considered in the subsections \ref{inhom1} and \ref{inhom2} below.

\subsection{Homogeneous equation with bounded boundary data} \label{inhom1}

We study the sub-solutions, super-solutions and solutions of $Lu=0$ in $U$ but $u$ needs not be zero on $\gamma$.


\begin{theorem}\label{thmih1}
The following statements hold true.
    \begin{enumerate}[label=\tnum]
        \item Assume $Lu\le 0$ in $U$. Then
        \beq\label{uMB0}
        \sup_{x\in \widetilde U} u(x)\le \sup_{x\in \mathcal S_{\theta_0}\cup\gamma} u(x),
        \eeq
        \beq\label{uMBa}
\max_{x\in \mathcal S_\theta}u(x)\le M_\gamma(1-C_* e^{-\nu(\theta-\theta_0)})+C_*\left(\sup_{x\in \mathcal S_{\theta_0}\cup\gamma} u(x)\right)e^{-\nu(\theta-\theta_0)}
\text{ for all $\theta\ge\theta_0$.}
\eeq
   Consequently,
        \beq\label{ublim1}
    \limsup_{\theta\to\infty}\left(\max_{x\in \mathcal S_\theta}u(x)\right)\le  \limsup_{\theta\to\infty}M_\gamma(\theta).
    \eeq

\item Assume $Lu\ge 0$ in $U$. Then
        \beq\label{umB0}
        \inf_{x\in \widetilde U} u(x)\ge \inf_{x\in \mathcal S_{\theta_0}\cup\gamma} u(x),
        \eeq
        \beq\label{uMBb}
\min_{x\in \mathcal S_\theta}u(x)\ge m_\gamma(1-C_*e^{-\nu(\theta-\theta_0)})+C_*\left(\inf_{x\in \mathcal S_{\theta_0}\cup\gamma} u(x)\right)e^{-\nu(\theta-\theta_0)}\text{ for all $\theta\ge\theta_0$.}
\eeq
   Consequently,
        \beqs
    \liminf_{\theta\to\infty}\left(\min_{x\in \mathcal S_\theta}u(x)\right)\ge  \liminf_{\theta\to\infty}m_\gamma(\theta).
    \eeqs
 \end{enumerate}
\end{theorem}
\begin{proof}
(i) Let $v=u-M_\gamma$. Then $Lv\le 0$ in $U$ and $v\le 0$ on $\gamma$.
Also $v$ is below exponential growth in $\widetilde U$.
Applying Corollary \ref{cor2}  to  $v$, we have Theorem \ref{thm1}\ref{T1i}\ref{T1ia}, which results in
\beqs
v(x)\le \max_{x\in \mathcal S_{\theta_0}}(u(x)-M_\gamma)^+ \text{ for any }x\in\widetilde U,
\eeqs
\beqs
v(x)\le C_*\left[\max_{x\in \mathcal S_{\theta_0}}(u(x)-M_\gamma)^+\right] e^{-\nu(\theta-\theta_0)} 
\text{ for any $\theta\ge\theta_0$ and $x\in \mathcal S_\theta$.}
\eeqs
These imply
\beq\label{uMall}
u(x)\le M_\gamma+\max_{x\in \mathcal S_{\theta_0}}(u(x)-M_\gamma)^+ \text{ for any }x\in\widetilde U,
\eeq
\beq\label{uMB}
u(x)\le M_\gamma+C_*\left[\max_{x\in \mathcal S_{\theta_0}}(u(x)-M_\gamma)^+\right]e^{-\nu(\theta-\theta_0)}
\text{ for any $\theta\ge\theta_0$ and $x\in \mathcal S_\theta$.}
\eeq
Observe that 
\begin{align*}
    \max_{x\in \mathcal S_{\theta_0}}(u(x)-M_\gamma)^+
    &=\max_{x\in \mathcal S_{\theta_0}} \Big(\max\{u(x),M_\gamma\}-M_\gamma\Big)
    =\max_{x\in \mathcal S_{\theta_0}} \Big(\max\{u(x),M_\gamma\}\Big)-M_\gamma\\
    &= \left(\max\left\{\max_{x\in \mathcal S_{\theta_0}} u(x),M_\gamma\right\}\right)-M_\gamma=\sup_{x\in \mathcal S_{\theta_0}\cup\gamma} u(x)-M_\gamma.
\end{align*}
Combining this with \eqref{uMall} and \eqref{uMB}, we obtain \eqref{uMB0} and \eqref{uMBa}.

Now, let 
$\widetilde M_\gamma(\theta)=\sup\{ u(x):x\in \gamma\cap \widetilde U_{[\theta,\infty)}\}$.
Given any number $\eta\ge \theta_0$. We apply inequality \eqref{uMBa} to $\theta\ge \eta$, with $\theta_0$ being replaced with $\eta$, and notice that $C_*$ and $\nu$ can stay the same. We obtain 
\beqs
\max_{x\in \mathcal S_\theta}u(x)\le \widetilde M_\gamma(\eta)(1-C_*e^{-\nu(\theta-\eta)})+C'e^{-\nu(\theta-\eta)},
\eeqs
where $C'$ depends on $u$ and $\eta$. Passing $\theta\to\infty$ gives 
\beqs
\limsup_{\theta\to\infty}\left(\max_{x\in \mathcal S_\theta}u(x)\right)\le \widetilde M_\gamma(\eta).
\eeqs
Then passing $\eta\to\infty$ gives
\beqs
\limsup_{\theta\to\infty} \left(\max_{x\in \mathcal S_\theta}u(x)\right)\le \lim_{\eta\to\infty}\widetilde M_\gamma(\eta)
\eeqs
which proves \eqref{ublim1}.

(ii) Applying part (i) to the function $(-u)$, we obtain part (ii). We omit the details. 
\end{proof}

From Theorem \ref{thmih1}, we obtain estimates for the oscillations next.

\begin{corollary}\label{corih1}
    Let $Lu=0$ in $U$. Then
    \beq\label{osc0}
     \osc_{x\in \widetilde U}u(x) \le \osc_{x\in \mathcal S_{\theta_0}\cup\gamma}u(x),
    \eeq
\beq\label{osc1}
    \osc_{x\in \mathcal S_\theta}u(x) \le\left(\osc_{x\in \gamma}u(x)\right)(1-C_*e^{-\nu(\theta-\theta_0)}) + C_*\left(\osc_{x\in \mathcal S_{\theta_0}\cup\gamma}u(x)\right) e^{-\nu(\theta-\theta_0)}\text{ for all $\theta\ge\theta_0$.}
    \eeq
    Consequently,
        \beq\label{osc2}
    \limsup_{\theta\to\infty}\left(\osc_{x\in \mathcal S_\theta}u(x)\right)\le     \limsup_{\theta\to\infty}M_\gamma(\theta)-\liminf_{\theta\to\infty}m_\gamma(\theta).
    \eeq
\end{corollary}
\begin{proof}
From \eqref{uMB0} and \eqref{umB0}, we have \eqref{osc0}.
From \eqref{uMBa} and \eqref{uMBb}, we have
\beq\label{oMmB}
 \osc_{x\in \mathcal S_\theta}u(x) \le (M_\gamma-m_\gamma)(1-C_*e^{-\nu(\theta-\theta_0)})+C_*\left(\osc_{x\in \mathcal S_{\theta_0}\cup\gamma} u(x)\right)e^{-\nu(\theta-\theta_0)}.
\eeq
Thus, we obtain \eqref{osc1}. Now, let 
 \beqs 
\widetilde M_\gamma(\theta)=\sup\{ u(x):x\in \gamma\cap \widetilde U_{[\theta,\infty)}\},\quad 
\widetilde m_\gamma(\theta)=\inf\{ u(x):x\in \gamma\cap \widetilde U_{[\theta,\infty)}\}.
\eeqs
Given $\eta\ge \theta_0$, applying inequality \eqref{oMmB} to $\theta\ge \eta$, with $\theta_0$ being replaced with $\eta$, yields
\beqs
\osc_{x\in \mathcal S_\theta}u(x)\le (\widetilde M_\gamma(\eta)-\widetilde m_\gamma(\eta))(1-C_*e^{-\nu(\theta-\eta)})+C'e^{-\nu(\theta-\eta)},
\eeqs
where $C'$ depends on $u$ and $\eta$. Passing $\theta\to\infty$ gives 
\beqs
\limsup_{\theta\to\infty}\left(\osc_{x\in \mathcal S_\theta}u(x)\right)\le \widetilde M_\gamma(\eta)-\widetilde m_\gamma(\eta).
\eeqs
Then passing $\eta\to\infty$ gives
\beqs
\limsup_{\theta\to\infty} \left(\osc_{x\in \mathcal S_\theta}u(x)\right)\le \lim_{\eta\to\infty}\widetilde M_\gamma(\eta)-\lim_{\eta\to\infty}\widetilde m_\gamma(\eta)
\eeqs
which proves \eqref{osc2}.
\end{proof}

\subsection{Inhomogeneous equation with a bounded forcing function} \label{inhom2}
Let $f:U\to\R$ be a given bounded function.
Define
$M_B=\sup_{x\in\gamma} |u(x)|$,
$M_f=\sup_{x\in U} |f(x)|$,  and
\beqs
M_B(\theta)=\max \left\{ |u(x)|:x\in\gamma\cap \mathcal S_\theta\right\},\quad 
M_f(\theta)=\sup \left\{ |f(x)|:x\in U\cap \mathcal S_\theta\right\}.
\eeqs

\begin{theorem}\label{thmih3}
    Suppose $Lu=f$ in $U$.  Then there exists a positive constant $K_1$ independent of $u$ and $f$ such that, for any $\theta\ge\theta_0$,
\beq\label{uMBF1}
\max_{x\in \mathcal S_\theta} |u(x)|\le M_B+K_1 M_f + C_*\left( \max_{x\in \mathcal S_{\theta_0}}|u(x)|+M_B+K_1 M_f \right)e^{-\nu(\theta-\theta_0)}.
\eeq
Consequently,
\beq\label{uMBF2}
\limsup_{\theta\to\infty}\left(\max_{x\in \mathcal S_\theta}|u(x)|\right)
\le \limsup_{\theta\to\infty}M_B(\theta)+K_1 \limsup_{\theta\to\infty}M_f(\theta).
\eeq 
\end{theorem}
\begin{proof}
Let $ d_1=\sup\{|x_1|:x=(x_1,x_2)\in U\}$ which is a positive number.
Set the positive numbers
\beq\label{KeK}
K_0=\frac{1+M_2}{c_0},\
\varep=e^{-K_0d_1}\text{ and  } 
K_1=\frac1{\varep e^{-K_0 d_1} K_0}=\frac{e^{2K_0d_1}}{K_0}.
\eeq 
For $x=(x_1,x_2)\in \widetilde U$, define
\beqs 
v(x)=u(x)-M_B-K_1M_f (1-\varep e^{K_0 x_1}) .
\eeqs 
 Denote the unit vector $e_1=(1,0)$. 
We have
\begin{align*}
Lv&=Lu-a_{11} K_1\varep M_f K_0^2 e^{K_0 x_1} + K_1 \varep M_f K_0 e^{K_0x_1} b(x)\cdot e_1  \\
&=f-\varep K_0K_1M_f( a_{11}K_0  - b(x)\cdot e_1) e^{K_0x_1}.
\end{align*}
By Assumption \ref{condB}, the definition of $M_f$, the fact $a_{11}(x)\ge c_0$ and the choice in \eqref{KeK}, we have
\begin{align*}
Lv&\le M_f-\varep K_1 K_0M_f( c_0 K_0  - M_2) e^{K_0x_1}
=M_f(1-\varep K_1 K_0 e^{K_0x_1})\\
&\le M_f(1-\varep K_1K_0 e^{-K_0d_1})=0.
\end{align*}
Thus, $Lv\le 0$ in $U$.
Note that 
$1-\varep e^{K_0 x_1}\ge 0$ in $\widetilde U$, hence, $v\le 0$ on $\gamma$.
Since solution $u$ is below exponential growth in $\widetilde U$, so is $v$.
We have from \eqref{d1a} of Theorem \ref{thm1} applied to $v$ that, for any $\theta\ge\theta_0$ and $x\in \mathcal S_\theta$,
\beqs
v(x)\le C_*\left(\max_{x\in \mathcal S_{\theta_0}}|v(x)| \right)e^{-\nu(\theta-\theta_0)},
\eeqs
thus,
\beq\label{uMBF}
u(x)\le M_B+K_1 M_f + C_*\left( \max_{x\in \mathcal S_{\theta_0}}|u(x)|+M_B+K_1 M_f \right)e^{-\nu(\theta-\theta_0)}.
\eeq
Applying \eqref{uMBF} to $u:=-u$ and $f:=-f$, and noticing that $M_B$ is the same for both $u$ and $-u$, and $M_{-f}=M_f$, we have, for any $\theta\ge\theta_0$ and $x\in \mathcal S_\theta$,
\beq\label{uMBFp}
-u(x)\le M_B+K_1 M_f + C_*\left( \max_{x\in \mathcal S_{\theta_0}}|u(x)|+M_B+K_1 M_f \right)e^{-\nu(\theta-\theta_0)}.
\eeq
Then the estimate \eqref{uMBF1} follows from \eqref{uMBF} and \eqref{uMBFp}.

We prove \eqref{uMBF2} next. Define, for $\theta>\theta_0$, 
\beqs
\widetilde M_B(\theta)=\sup \left\{ |u(x)|:x\in\gamma\cap \widetilde U_{[\theta,\infty)}\right\},\quad 
\widetilde  M_f(\theta)=\sup \left\{ |f(x)|:x\in U_{[\theta,\infty)}\right\}.
\eeqs
Let $\eta>\theta_0$.  Applying inequality \eqref{uMBF1} to $\theta\ge \eta$, with $\eta$ replacing $\theta_0$, gives
\beqs
\max_{x\in \mathcal S_\theta} |u(x)|\le \widetilde M_B(\eta)+K_1 \widetilde M_f(\eta)+Ce^{-\nu(\theta-\eta)},
\eeqs
where $C>0$ depends on $u$, $\widetilde M_B(\eta)$ and $\widetilde M_f(\eta)$. Passing $\theta\to\infty$ gives 
\beqs
\limsup_{\theta\to\infty} \left(\max_{x\in \mathcal S_\theta} |u(x)|\right)\le \widetilde M_B(\eta)+K_1 \widetilde M_f(\eta).
\eeqs
Then passing $\eta\to\infty$ gives
\beqs
\limsup_{\theta\to\infty}\left(\max_{x\in \mathcal S_\theta} |u(x)|\right)\le \lim_{\eta\to\infty}\widetilde M_B(\eta)+K_1 \lim_{\eta\to\infty}\widetilde  M_f(\eta)
\eeqs
which proves \eqref{uMBF2}.
\end{proof}

\section{Uniqueness and continuous dependence}\label{unicd}

In this section, let Assumption \ref{firstA} hold for $\Omega=U$.
Let $u_1$ and $u_2$ be functions in $C^2(U)\cap C(\widetilde U)$.
Recall that $C_*$ and $\nu$ are the positive constants in Theorem \ref{thm1},
 $\Theta_0\ge \bar\Theta$ and $\bar\nu>0$ are the numbers in Theorem \ref{thm3}, and $K_1>0$ is the number in Theorem \ref{thmih3}.

\begin{theorem}\label{ThmDep1}
Suppose $Lu_1=Lu_2$ in $U$ and $u_1=u_2$ on $\gamma$. 
\begin{enumerate}[label=\tnum]
    \item Under Assumption \ref{condB} for $\Omega=U$, if   $u_1$, $u_2$ are below exponential growth in $\widetilde U$, then
    \beq\label{dineq1}
\sup_{x\in U} |u_1(x)-u_2(x)|\le \max_{x\in \mathcal S_{\theta_0}}|u_1(x)-u_2(x)|,
\eeq
\beq\label{dineq2}
\max_{x\in \mathcal S_{\theta}} |u_1(x)-u_2(x)|\le C_*\left(\max_{x\in \mathcal S_{\theta_0}}|u_1(x)-u_2(x)| \right)e^{-\nu(\theta-\theta_0)}
\text{ for all $\theta\ge\theta_0$.}
\eeq

\item Under Assumption \ref{bmassum}, if $u_1$ and $u_2$ are bounded in $\widetilde U$, then
\beq\label{dineq3}
\max_{x\in \mathcal S_{\theta}}|u_1(x)-u_2(x)|\le \left(\max_{x\in \mathcal S_{\theta_0}}|u_1(x)-u_2(x)|\right) \min\left\{1,\exp\left( -\bar\nu\int_{\Theta_0}^\theta \kappa^{m(\zeta)} \d\zeta\right)\right\}.
\eeq
\end{enumerate}
\end{theorem}
\begin{proof}
    Let $v=u_1-u_2$, then $Lv=0$ in $U$ and $v=0$ on $\gamma$.
    We  apply Theorem \ref{thm1} \ref{T1ii}, \ref{T1iia}  and Theorem \ref{thm4} \ref{T4i} to the solution $v$. Then \eqref{dineq1}, \eqref{dineq2} and \eqref{dineq3} follow from \eqref{dabs0},  \eqref{dabsa} and \eqref{ubdd1}, respectively.
\end{proof}

In the next theorem, we show an additional dependence  on the forcing function.

\begin{theorem}\label{ThmDep2}
Assume
\begin{itemize}
    \item $b(x)$ satisfies Assumption \ref{condB} for $\Omega=U$, and two given functions  $f_1,f_2:U\to \R$ are bounded in $U$,
    \item $Lu_1=f_1$ and $Lu_2=f_2$ in $U$,  
    \item $u_1$, $u_2$ are below exponential growth in $\widetilde U$ and are bounded on $\gamma$. 
\end{itemize}
Then one has   
\beq\label{udep3}
\begin{aligned}
 \max_{x\in \mathcal S_{\theta}}|u_1(x)-u_2(x)|
&\le (1+C_*)\sup_{x\in \gamma}|u_1(x)-u_2(x)| +K_1(1+C_*) \sup_{x\in U}|f_1(x)-f_2(x)|\\
&\quad
+ C_*\left( \max_{x\in \mathcal S_{\theta_0}}|u_1(x)-u_2(x)|\right)e^{-\nu(\theta-\theta_0)}
\text{ for all $\theta\ge\theta_0$,}
\end{aligned}
\eeq    
\beq\label{udep4}
\begin{aligned}
\limsup_{\theta\to\infty}\left(\max_{x\in \mathcal S_{\theta}}|u_1(x)-u_2(x)|\right)
&\le \limsup_{\theta\to\infty}\left(\sup_{x\in \gamma\cap \mathcal S_\theta}|u_1(x)-u_2(x)|\right)\\
&\quad +K_1 \limsup_{\theta\to\infty}\left(\sup_{x\in U\cap \mathcal S_\theta}|f_1(x)-f_2(x)|\right).
\end{aligned}
\eeq 
\end{theorem}
\begin{proof}
    Let $v=u_1-u_2$ and $f=f_1-f_2$. Then $Lv=f$ in $U$ and $v$ is below exponential growth in $\widetilde U$. We apply Theorem \ref{thmih3} to the solution $v$ and obtain \eqref{udep3} and  \eqref{udep4}  from \eqref{uMBF1} and  \eqref{uMBF2}, respectively.
\end{proof}

\begin{remark}\label{adjust} The results obtained in this paper can be easily adapted to the following slightly different situations.
(a) The domain can wind around the origin clockwise instead of counter-clockwise as in section \ref{domsec}.
(b) Through winding, the domain can approach the circle $\mathcal C_*$ from the inside instead of from the outside.
\end{remark}

\begin{remark}\label{more}
Our study above is focused on the qualitative properties of the solutions without discussing their existence. For that matter, see related results, for example,  in \cite{LaNa1985} for positive solutions in a conic domain, or \cite{IbraPhd1985} for unbounded solutions in unbounded domains.
This research can also be expanded into many directions. 
\begin{enumerate}[label=\rnum]
    \item The nonlinear problem such as the elliptic version of the models in \cite{HI3,HI4}. As seen in the cited papers, the analysis of the linear problems like the one in the current work is, in fact, very crucial.
    \item The domain can spiral around a more general limiting set such as a general loop or even a point.
    \item The domain can have more complicated geometry but lies inside our domain $U$ above.
\end{enumerate}
\end{remark}

\appendix
\section{}\label{apx}

\begin{proof}[Proof of Lemma \ref{subsol}]
Set $u(x)=\varphi(|x-x_*|)$, where $\varphi$ is  a function in the class $C^2((0,\infty),\R)$ with $\varphi'\le 0$ in $(0,\infty)$. 
  Clearly, $u\in C^2(\Omega)$. For $x\in \Omega$, denote $r=r(x)\eqdef |x-x_*|$ and $\xi=\xi(x)\eqdef (x-x_*)/r$.
Denote $a_*(x)=\sum_{i,j=1}^n a_{ij}(x)\xi_i(x)\xi_j(x)=\xi^{\rm T}(x)A(x)\xi(x)$.  
Note that $a_*(x)$ is positive thanks to \eqref{welip}.
Elementary calculations give 
\begin{align*}
Lu(x)
&=-\varphi''(r)\sum_{i,j=1}^n a_{ij}(x)\xi_i\xi_j-\frac{\varphi'(r)}{r}\sum_{i,j=1}^n a_{i,j}(x)(\delta_{ij}-\xi_i\xi_j)+\frac{\varphi'(r)}{r} b(x)\cdot (x-x_*)\\
&= -  a_*(x)\left(\varphi''(r)+\frac{\varphi'(r)}{r}\left[\frac{{\rm Tr}A(x)- b(x)\cdot (x-x_*)}{a_*(x)}-1\right]\right). 
\end{align*}
By \eqref{b-cond} and the facts $-a_*(x)\varphi'(r)\ge 0$, it follows that
\beqs
Lu(x)\le -a_*(x)\left (\varphi''(r)+\frac{\varphi'(r)}{r}(e_0-1 )\right).
\eeqs
Therefore, a sufficient condition for $Lu(x)\le 0$ is 
$\varphi''(r)+r^{-1}\varphi'(r)(e_0-1 )\ge 0$.
Now, by choosing $\varphi(t)=t^{-s}$, the last condition becomes $r^{-s-2}s(s+2-e_0)\ge 0.$
Clearly, this requirement is satisfied when $s>0$ and  $s\ge e_0-2$.
\end{proof}

\begin{lemma}[Dichotomy for sequences]\label{maindiseq}
    Let $(a_i)_{i=0}^\infty$ be a sequence in $[0,\infty)$ and $(\lambda_i)_{i=1}^\infty$ be a sequence in $(0,1)$  such that 
    \beq\label{seq1}
    a_i\le \lambda_i\max\{a_{i-1},a_{i+1}\} \text{ for all } i\ge 1.
    \eeq
Then there are only the following two possibilities.
\begin{enumerate}[label=\tnum]
 \item One has 
 \beq\label{adec} 
 a_i\le \lambda_i a_{i-1} \text{ for all }i\ge 1.
 \eeq 
 Consequently, for all $i\ge 1$,
 \beq\label{seq6}
 a_{i}\le \lambda_{1}\lambda_{2}\ldots \lambda_{i}a_0.
 \eeq  
 \item There exists $i_*\ge 1$ such that 
 \beq \label{seq7}
a_{i_*}>0 \text{ and } a_{i}\ge \frac{a_{i-1}}{\lambda_{i-1}} 
  \text{ for all $i\ge i_*+1$.}
 \eeq  
Consequently,
 \beq \label{seq5}
a_{i}\ge \frac{a_{i_*}}{\lambda_{i_*}\lambda_{i_*+1}\ldots \lambda_{i-1}} 
  \text{ for all $i\ge i_*+1$.}
 \eeq  
 \end{enumerate}
\end{lemma}
\begin{proof}
We refine the argument in {\cite[Lemma V.6]{HIK2}}.
First of all, note from \eqref{seq1}, for any $i\ge 1$, that we have either of the following cases.
 \begin{enumerate}[label=\rnum]
  \item\label{casa} $a_{i-1}< a_{i+1}$ which, together with \eqref{seq1}, implies $a_{i}\le \lambda_ia_{i+1}$. 
  \item\label{casb} $a_{i-1}\ge a_{i+1}$  which, together with \eqref{seq1}, implies $a_{i}\le \lambda_ia_{i-1}$.
 \end{enumerate}
We have only the following two cases to consider.

 \emph{Case 1. For all $i\ge 1$,   Case \ref{casb} holds.} Then we have 
 $ a_{i} \le \lambda_{i}a_{i-1}$ for all $i\ge 1$ which proves \eqref{adec}.
Applying \eqref{adec} repeatedly, one has, for any $i\ge 1$,
    \beqs
a_{i} \le \lambda_{i} a_{i-1} \le \lambda_{i}\lambda_{i-1} a_{i-2} 
\le \ldots\le \lambda_{i}\lambda_{i-1} \ldots \lambda_{1}a_{0}
  \eeqs 
which proves \eqref{seq6}.

\emph{Case 2. There is $i_0\ge 1$ such that Case \ref{casa} holds for $i=i_0$.} Then we have
 \beq\label{a}  
  a_{i_0-1}< a_{i_0+1} \quad \text { and } \quad a_{i_0}\le \lambda_{i_0}a_{i_0+1}.
 \eeq
Consider  $i=i_0+1$ and suppose Case \ref{casb} holds. Then we have  
\beqs
a_{i_0} \ge   a_{i_0+2}\text{  and  }a_{i_0+1}\le \lambda_{i_0+1}a_{i_0}.
\eeqs 
Together with \eqref{a}, the last inequality yields   
\beq\label{aii}
a_{i_0+1}\le \lambda_{i_0+1}a_{i_0} \le \lambda_{i_0+1}\lambda_{i_0} a_{i_0+1}.
\eeq 
Since $a_{i_0+1}\ge 0$ and $\lambda_{i_0+1}\lambda_{i_0} <1$, inequality \eqref{aii} implies $a_{i_0+1}=0$. 
This fact and the first inequality in \eqref{a} imply $a_{i_0-1}<0$ which is a contradiction.
Therefore, we must have Case \ref{casa} for $i=i_0+1$.
Then, by induction, Case \ref{casa} holds for $i=i_0+j$ for all $j\ge 0$. Consequently,
\beq\label{aincr}
a_{i}\le \lambda_{i}a_{i+1}\text{ for all }i\ge i_0.
\eeq

Set $i_*=i_0+1$. Note from \eqref{a} that $a_{i_*}>a_{i_0-1}\ge 0$, hence we obtain the first part of \eqref{seq7}. The second part of  \eqref{seq7} follows directly from \eqref{aincr}.
For any $i\ge i_*+1$, applying \eqref{aincr} repeatedly yields
\beqs
a_{i_*}\le \lambda_{i_*}a_{i_*+1}\le \lambda_{i_*}\lambda_{i_*+1}a_{i_*+2}
\le\ldots\le \lambda_{i_*}\lambda_{i_*+1}\ldots \lambda_{i-1} a_{i},
\eeqs 
which implies \eqref{seq5}.   
\end{proof}

\medskip
\noindent\textit{Acknowledgment.} A.I. obligatorily acknowledges OGRI's grant 122022800272-4. 


\medskip
\noindent\textbf{Data availability.} 
No new data were created or analyzed in this study.


\medskip
\noindent\textbf{Conflict of interest.}
There are no conflicts of interests.

\bibliography{paperbaseall}{}
\bibliographystyle{plain}
\end{document}